\documentclass[11pt]{article}
\usepackage{amssymb}
\usepackage{amsmath}
\usepackage{amsthm}
\usepackage{latexsym}
\usepackage{graphicx}
\usepackage{enumerate}

\numberwithin{equation}{section}

\def\rc#1{\frac{1}{#1}}
\def\qed{\square}
\def\cases#1{\left\{\begin{matrix} #1 \end{matrix}\right.}
\def\pmatrix#1{\left(\begin{matrix} #1 \end{matrix}\right)}
\def\nonpmatrix#1{\begin{matrix} #1 \end{matrix}}

\newtheorem{Lemma}{Lemma}
\newtheorem{Corollary}{Corollary}

\newtheorem{Theorem}{Theorem}
\newtheorem*{remark}{Remark}
\newcommand{\ba}[1]{\begin{array}{@{}#1@{}}}
\newcommand{\ea}{\end{array}}
\long\def\symbolfootnote[#1]#2{\begingroup%
\def\thefootnote{\fnsymbol{footnote}}\footnote[#1]{#2}\endgroup}

\begin{document}
\newcommand{\ul}{\underline}
\newcommand{\be}{\begin{equation}}
\newcommand{\ee}{\end{equation}}
\newcommand{\ben}{\begin{enumerate}}
\newcommand{\een}{\end{enumerate}}

\long\def\symbolfootnote[#1]#2{\begingroup%
\def\thefootnote{\fnsymbol{footnote}}\footnote[#1]{#2}\endgroup}
\def\rc#1{\frac{1}{#1}}
\def\qed{\square}
\def\cases#1{\left\{\begin{matrix} #1 \end{matrix}\right.}
\def\pmatrix#1{\left(\begin{matrix} #1 \end{matrix}\right)}
\def\nonpmatrix#1{\begin{matrix} #1 \end{matrix}}
\def\emx{{\cal E}}
\def\jmx{{\cal J}}
\def\imx{{\cal I}}
\def\omx{{\cal O}}
\def\xmx{{\cal X}}
\def\umx{{\cal U}}
\par

\title{On Block Representations and Spectral Properties of %Magic and
 Semimagic Square Matrices% and Quadratic Forms
}
\date{\today}

\author{S. L. Hill, M. C. Lettington\\ and K. M. Schmidt (Cardiff)$$}
\maketitle
\begin{abstract}
Using the decomposition of semimagic squares into the associated and
balanced symmetry types as a motivation, we introduce an equivalent representation
in terms of block-structured matrices.
This block representation provides a way of constructing such matrices with
further symmetries and of studying their algebraic behaviour, significantly
advancing and contributing to the understanding of these symmetry properties.
In addition to studying classical attributes, such as dihedral equivalence
and the spectral properties of these matrices, we show that the inherent
structure of the block representation facilitates the definition of low-rank semimagic
square matrices. This is achieved by means of tensor product blocks. Furthermore, we study the rank and
eigenvector decomposition of these matrices, enabling the construction of a corresponding
two-sided eigenvector matrix in rational terms of their entries. The paper concludes with the
derivation of a correspondence between the tensor product block representations and
quadratic form expressions of Gaussian type.
\end{abstract}

\section{Introduction}
\noindent
The magic square of order 3 or ``Loh Shu square'' \cite{han} was known in China as early as the Warring States period (481 BC -- 221 BC).
\begin{table}[h]
\be
\begin{array}{|c|c|c|}
\hline {}
8 & 1 & 6 \\
\hline
3 & 5 & 7 \\
\hline
 4 & 9 & 2 \\
 \hline
\end{array}
%\label{fig:m8.1}
\nonumber
\ee
%\begin{center}
%\text{ The unique (up to dihedral symmetry) $3\times 3$ square.}
%\end{center}
\end{table}
Up to rotations and reflections, it gives the unique arrangement of the numbers
$1, \dots, 9$ in a $3 \times 3$ square such that all the rows and columns and both
diagonals add up to the same number (in this case, $15$).
More generally, an arrangement of any $n^2$ numbers in an $n \times n$ square such
that all rows, columns and both diagonals add up to the same number is called a
magic square; a semimagic square if the requirement on the diagonals is dropped
\cite{andrewsmagic}, \cite{essen}.

A.~C.~Thompson \cite{oldmagic} was one of the first to publish on the subject of matrix multiplication of magic squares, observing that the above $3\times 3$ magic square, considered as a square matrix, retains its symmetry for all odd integer powers.
We know now that this property is due to the fact that this is a so-called associated
semimagic square (or type A matrix) and that the type A matrices form a
$\mathbb{Z}_2$-graded algebra with the so-called balanced semimagic squares (or
type B matrices) of the same dimension;
such that a product of any two type A matrices gives a type B matrix, while the
product of a type A and a type B matrix is of type A, and the type B matrices form
a matrix subalgebra by themselves (cf. \cite{mcl1}, \cite{rBLS} Lemma 1.1).
Any semimagic square matrix can be written as a sum of an associated and a balanced
semimagic square matrix, so the above graded algebra is in fact the matrix algebra of
semimagic squares (\cite{rBLS} Lemma 2.1).

In this paper we introduce an equivalent {\it block representation\/} of $n\times n$
semimagic squares by means of conjugation with a symmetric involution matrix $\xmx_n$;
in the transformed representation, split into a $2 \times 2$ array of block matrix
components, the types A and B correspond to the off-diagonal and on-diagonal blocks,
respectively, so the symmetry type decomposition of the semimagic square matrix can
easily be read off. Moreover, it turns out that the block representation reveals
other traits of the semimagic square and makes them more accessible to analysis than
the original form. For example, based on the groundwork presented here, the authors
have studied a wider range of matrix symmetries, including the type of `most perfect pandiagonal magic squares' \cite{pandiag}, using the block representation,
identifying further structured algebras of semimagic square (and related) matrices
\cite{algebra}.

The present work is structured as follows.
In Section 2, we introduce the block representation and derive its properties
equivalent to the symmetry types A and B; there is a curious fundamental difference
between the cases of even-dimensional and odd-dimensional matrices. As
the block representation is an algebra isomorphism which is covariant under transposition, it is possible to
study algebraic properties by looking at the transformed matrices instead of the
original semimagic squares. Moreover, the dihedral operations of rotation and
reflection, typically resulting in different matrices associated with essentially
the same semimagic square, are shown to be expressed by a simple sign change in
different component blocks.

In Section 3, we address the question of how magic squares can be recognised in
terms of their block representation; as the type A part is always magic, this is
a question for type B matrices only. As a result, we show that for any balanced
magic square matrix, some power of it will be either trivial (i.e. with all matrix
entries the same) or non-magic.

In Section 4, we consider the existence of inverses or quasi-inverses of
semimagic square matrices, and more specifically consider the ranks and
eigenvector decomposition of associated magic square matrices and their
(balanced) matrix algebra squares. Here the concept of {\it parasymmetry\/}, i.e.
the property of a matrix to have a symmetric square, plays a crucial role.
The case of matrices whose block representation is composed of rank 1 blocks
allows a particularly detailed analysis, including the explicit construction of
a two-sided (both right and left) eigenvector matrix.

In the final section, we use the latter results to establish a connection to
quadratic forms of the type studied in number theory \cite{gauss}.

\section{Block Representations and Dihedral Symmetries}
{\bf Definition.}
An
$n \times n$
matrix
$M \in {\mathbb R}^{n \times n}$
is called
{\it semimagic square of weight\/}
$w$
if each of its rows and columns sums to
$n w.$
If, in addition, each of its two diagonals also adds up to
$n w,$
it is called a
{\it magic square.\/}
\par\medskip\noindent
The following two centre-point symmetry types of semimagic square matrices are of interest.
\par\medskip\noindent
{\bf Definition.}
Let
$M$
be an
$n \times n$
semimagic square matrix of weight
$w.$
\begin{description}
\item{(a)}
The matrix
$M$
is called
{\it associated\/}
if each entry and its mirror entry w.r.t. the centre of the matrix add to
$2w,$
i.e. if
\begin{align}
 M_{i j} + M_{n+1-i, n+1-j} &= 2 w
 \qquad (i,j \in \{1, \dots, n\}).
\nonumber\end{align}
\item{(b)}
The matrix
$M$
is called
{\it balanced\/}
if each entry is equal to its mirror entry w.r.t. the centre of the matrix, i.e.
\begin{align}
 M_{i j} - M_{n+1-i, n+1-j} &= 0
 \qquad (i,j \in \{1, \dots, n\}).
\nonumber\end{align}
\end{description}
\par\medskip\noindent
These properties can be equivalently expressed in a more convenient way using the
following conventions.
Let
$1_n \in {\mathbb R}^n$
be the vector which has all entries equal to
$1.$
Let
$\emx_n = (1)_{i,j = 1}^n \in {\mathbb R}^{n \times n}$
be the
$n \times n$
matrix which has all entries equal to
$1.$
Moreover, let
$\jmx_n = (\delta_{i, n+1-j})_{i,j = 1}^n \in {\mathbb R}^{n \times n},$
where
$\delta$
is the Kronecker symbol,
be the
$n \times n$
matrix which has entries
$1$
on the antidiagonal and
$0$
otherwise.
Later, we shall also use the notations
$0_n$
for the null vector in
${\mathbb R}^n,$
$\omx_n = (0)_{i, j = 1}^n$
for the
$n \times n$
null matrix, and
$\imx_n = (\delta_{i j})_{i,j = 1}^n$
for the
$n \times n$
unit matrix.
\par
Note that we consider the elements of
${\mathbb R}^n$
as column vectors, i.e. matrices with
$n$
rows
and a single column, throughout; row vectors are represented in the form
$v^T,$
where
$v \in {\mathbb R}^n.$
\begin{Lemma}\label{lsemi}
Let
$M \in {\mathbb R}^{n \times n}.$
\begin{description}
\item{(a)}
The following statements are equivalent.
\begin{description}
\item{(i)}
$M$
is semimagic square of weight
$w;$
\item{(ii)}
$M \emx_n = n w \emx_n = M^T \emx_n;$
\item{(iii)}
$1_n$
is an eigenvector, with eigenvalue
$n w,$
for both
$M$
and its transpose
$M^T.$
\end{description}
\item{(b)}
If
$M$
is semimagic square, then it is associated if and only if
$M + \jmx_n M \jmx_n = 2 w \emx_n.$
%\phantom{.}\hfill$(@@\egacond)$\par\noindent
\item{(c)}
If
$M$
is semimagic square, then it is balanced if and only if
$M = \jmx_n M \jmx_n.$
%\phantom{.}\hfill$(@@\egbcond)$\par\noindent
\end{description}
\end{Lemma}

\par\medskip\noindent
Let
$S_n$
be the set of all
$n \times n$
semimagic square matrices, and
$A_n, B_n$
the subsets of associated and balanced semimagic square matrices, respectively.
Then
$S_n$
is a subalgebra of the standard
$n \times n$
matrix algebra, and
$A_n$
and
$B_n$
are vector subspaces of
$S_n.$
Moreover,
$B_n$
is a subalgebra of
$S_n,$
and
\begin{equation}
\label{ealg}
A_n A_n \subset B_n, \qquad
A_n B_n \subset A_n, \qquad
B_n A_n \subset A_n.
\end{equation}
Furthermore,
$A_n$
and
$B_n$
generate the whole algebra
$S_n$
in the following way.
Denoting by
$S_n^o$
the subset (in fact, subalgebra) of weight 0 semimagic square matrices, and setting
$A_n^o = A_n \cap S_n^o,$
$B_n^o = B_n \cap S_n^o,$
we have
\begin{align}
 A_n \cap B_n &= {\mathbb R} \emx_n,
\nonumber\\
 A_n^o \cap B_n^o &= \{\omx_n\},
\nonumber\\
 \hbox{\rm and} \quad S_n &= A_n^o + B_n^o + {\mathbb R} \emx_n,
\nonumber\end{align}
where every
$n \times n$
semimagic square matrix of weight
$w$
can be written as a sum of unique elements of
$A_n^o$
and
$B_n^o,$
and
$w \emx_n$
(cf. \cite{rBLS} Lemma 2.3).
%\par\bigskip\noindent{\bf {[sBlr]}. Block representation of semimagic square matrices}\nobreak
\par\medskip\noindent
\subsection*{Even Dimensional Case}
\par\medskip\noindent
We first consider
$2n \times 2n $
matrices,
$n \in {\mathbb N}.$
Let
\begin{align}
 \xmx_{2n} = \rc{\sqrt 2} \pmatrix{ \imx_n & \jmx_n \cr \jmx_n & -\imx_n} \in {\mathbb R}^{2n \times 2n}.
\nonumber\end{align}
Clearly
$\xmx_{2n}^2 = \imx_{2n}$
and
$\xmx_{2n}^T = \xmx_{2n},$
so
$\xmx_{2n}$
is an orthogonal symmetric involution.
Conjugation with the matrix
$\xmx_{2n}$
gives rise to a block representation of matrices in
$S_{2n}$
in which the symmetry type can easily be read off. This also provides a convenient and
systematic way of constructing semimagic square matrices with (or without) centre-point
symmetries.
\begin{Theorem}\label{tblre}
Let
$M \in {\mathbb R}^{2n \times 2n}.$
\begin{description}
\item{(a)}
The matrix
$M \in A_{2n}^o$
if and only if
\begin{align}
 M &= \xmx_{2n} \pmatrix{\omx_n & V^T \cr W & \omx_n} \xmx_{2n}
\nonumber\end{align}
where
$V, W \in {\mathbb R}^{n \times n}$
have row sum 0, i.e.
\begin{align}
 V 1_n &= 0_n, \qquad W 1_n = 0_n.
\nonumber\end{align}
\item{(b)}
The matrix
$M \in B_{2n}^o$
if and only if
\begin{align}
 M &= \xmx_{2n} \pmatrix{Y & \omx_n \cr \omx_n & Z} \xmx_{2n}
\nonumber\end{align}
where
$Y \in S_n^o$
and
$Z \in {\mathbb R}^{n \times n}.$
\item{(c)}
The weight matrix
$\emx_{2n}$
satisfies
\begin{align}
 \emx_{2n} &= \xmx_{2n} \pmatrix{2 \emx_n & \omx_n \cr \omx_n & \omx_n} \xmx_{2n}.
\nonumber\end{align}
\end{description}
\end{Theorem}
\par\medskip\noindent
In view of the decomposition of general semimagic square matrices mentioned above, we
see that, as a consequence of Theorem \ref{tblre}, any even-dimensional, weight
$w$
semimagic square matrix has, after
conjugation with
$\xmx_{2n},$
the block representation
\begin{align}
 \pmatrix{Y + 2 w \emx_n & V^T \cr W & Z},
\nonumber\end{align}
with a weight 0 semimagic square matrix
$Y,$
matrices
$V, W$
with row sum 0 and a matrix
$Z$
which can be any
$n \times n$
matrix. Evidently, this block representation clearly shows the decomposition into an
associated and a balanced matrix, corresponding to setting the two diagonal or the two
off-diagonal blocks equal to 0, respectively.
\par
From the block representation, it is very straightforward to generate all matrices in
$A_{2n}^o$
(and hence, by adding a multiple of
$\emx_{2n},$
in
$A_{2n});$
indeed, the conditions on the matrices
$V, W$
can very easily be satisfied, as
$n-1$
columns can be arbitrary when the last column is chosen so that the rows add to 0.
The construction of a general matrix in
$B_{2n}$
is a bit more complicated, as
$Z$
can be chosen arbitrarily, but
$Y$
must be a semimagic square matrix. At least this reduces the dimension of the problem from
$2n$
to
$n.$
\par
As the matrix
$\xmx_{2n}$
is symmetric (i.e., equal to its transpose) and involutory (i.e., its own inverse matrix),
the following observation is straightforward.
\begin{Lemma}\label{liso}
\begin{description}
\item{(a)}
The block representation is an algebra isomorphism; indeed,
\begin{align}
 \xmx_{2n} (\alpha M + N) \xmx_{2n} &= \alpha \xmx_{2n} M \xmx_{2n} + \xmx_{2n} N \xmx_{2n}
 \qquad (\alpha \in {\mathbb R}; M, N \in {\mathbb R}^{2n \times 2n})
\nonumber\end{align}
and
\begin{align}
 \xmx_{2n} (M N) \xmx_{2n} &= \xmx_{2n} M \xmx_{2n} \xmx_{2n} N \xmx_{2n}
 \qquad (M, N \in {\mathbb R}^{2n \times 2n}).
\nonumber\end{align}
\item{(b)}
The block representation of the transposed matrix is the transpose of the block
representation of the original matrix; indeed,
\begin{align}
 \xmx_{2n} M^T \xmx_{2n} &= (\xmx_{2n} M \xmx_{2n})^T
 \qquad (M \in {\mathbb R}^{2n \times 2n}).
\nonumber\end{align}
\end{description}
\end{Lemma}
\par\medskip\noindent
{\it Proof\/}
of Theorem \ref{tblre}.
We begin by writing the matrix
$M$
in the form
\begin{align}
 M &= \pmatrix{A & C \cr B & D},
\nonumber\end{align}
where
$A, B, C, D \in {\mathbb R}^{n \times n};$
then
\begin{align}
 \jmx_{2n} M \jmx_{2n} &= \pmatrix{\jmx_n D \jmx_n & \jmx_n B \jmx_n \cr \jmx_n C \jmx_n & \jmx_n A \jmx_n}.
\nonumber\end{align}
\begin{description}
\item{(a)}
Assume that
$M \in A_{2n}^o.$
Then, by Lemma \ref{lsemi} (b),
\begin{align}
 \omx_{2n} &= M + \jmx_{2n} M \jmx_{2n} = \pmatrix{\jmx_n (D + \jmx_n A \jmx_n) \jmx_n &
   \jmx_n (B + \jmx_n C \jmx_n) \jmx_n \cr B + \jmx_n C \jmx_n & D + \jmx_n A \jmx_n}
\nonumber\end{align}
(where we used
$\jmx_n^2 = \imx_n).$
This is equivalent to
$B = - \jmx_n C \jmx_n,$
$D = - \jmx_n A \jmx_n.$
Hence,
\begin{align}
 \xmx_{2n} M \xmx_{2n} &= \rc 2 \pmatrix{\imx_n & \jmx_n \cr \jmx_n & -\imx_n}
   \pmatrix{A & C \cr - \jmx_n C \jmx_n & -\jmx_n A \jmx_n}
   \pmatrix{\imx_n & \jmx_n \cr \jmx_n & -\imx_n}
\nonumber\\
 &= \pmatrix{ \omx_n & A \jmx_n - C \cr \jmx_n A + \jmx_n C \jmx_n & \omx_n}.
\label{eaeform}\end{align}
By Lemma \ref{lsemi} (a),
$M$
is a weight 0 semimagic square matrix if and only if
$1_{2n}$
is an eigenvector of both
$M$
and
$M^T$
for eigenvalue 0. In view of Lemma \ref{liso}, this is equivalent to
\begin{align}
 \xmx_{2n} 1_{2n} &= \sqrt 2 \pmatrix{1_n \cr 0_n}
\label{eonerep}\end{align}
being an eigenvector, for eigenvalue 0, of both
$\xmx_{2n} M \xmx_{2n}$
and
$(\xmx_{2n} M \xmx_{2n})^T.$
From (\ref{eaeform}), this corresponds to the conditions
\begin{align}
 \jmx_n (A + C \jmx_n) 1_n &= 0_n, \qquad (A \jmx_n - C)^T 1_n = 0_n.
\nonumber\end{align}
Thus
$V = (A \jmx_n - C)^T$
and
$W = \jmx_n A + \jmx_n C \jmx_n$
will have row sum 0.
\hfill\break
Conversely, given
$V, W \in {\mathbb R}^{n \times n}$
with row sum 0, we take
$A = \rc 2 (V^T \jmx_n + \jmx_n W),$
$C = \rc 2 (\jmx_n W \jmx_n - V^T)$
and further
$B = -\jmx_n C \jmx_n,$
$D = -\jmx_n A \jmx_n,$
to construct a weight 0 semimagic square matrix.
\item{(b)}
Assume
$M \in B_n^o.$
Then, by Lemma \ref{lsemi} (c),
\begin{align}
 \omx_{2n} &= M - \jmx_{2n} M \jmx_{2n} = \pmatrix{-\jmx_n (D - \jmx_n A \jmx_n) \jmx_n &
   -\jmx_n (B - \jmx_n C \jmx_n) \jmx_n \cr B - \jmx_n C \jmx_n & D - \jmx_n A \jmx_n},
\nonumber\end{align}
which is equivalent to
$B = \jmx_n C \jmx_n,$
$D = \jmx_n A \jmx_n.$
Hence
\begin{align}
 \xmx_{2n} M \xmx_{2n} &= \rc 2 \pmatrix{\imx_n & \jmx_n \cr \jmx_n & -\imx_n}
 \pmatrix{A & C \cr \jmx_n C \jmx_n & \jmx_n A \jmx_n}
 \pmatrix{\imx_n & \jmx_n \cr \jmx_n & -\imx_n}
\nonumber\\
 &= \pmatrix{A + C \jmx_n & \omx_n \cr \omx_n & \jmx_n A \jmx_n - \jmx C}.
\label{ebeform}\end{align}
As before, the condition that
$M$
is a semimagic square matrix of weight 0 means that
%$(@\eonerep)$
is an eigenvector with eigenvalue 0 of both
$\xmx_{2n} M \xmx_{2n}$
and
$(\xmx_{2n} M \xmx_{2n})^T;$
by (\ref{ebeform}) this is equivalent to
\begin{align}
 (A + C \jmx_n) 1_n &= 0_n, \qquad (A + C \jmx)^T 1_n = 0_n.
\nonumber\end{align}
By Lemma \ref{lsemi} (a),
$Y = A + C \jmx_n \in S_n^o.$
\hfill\break
Conversely, given
$Y \in S_n^o$
and
$Z \in {\mathbb R}^{n \times n},$
we take
$A = \rc 2 (Y + \jmx_n Z \jmx_n),$
$C = \rc 2 (Y \jmx_n - \jmx_n Z),$
and further
$B = \jmx_n C \jmx_n,$
$D = \jmx_n A \jmx_n,$
to construct a weight 0 balanced semimagic square matrix.
\item{(c)}
By a straightforward calculation,
\begin{align}
 \xmx_{2n} \emx_{2n} \xmx_{2n} &= \rc 2 \pmatrix {\imx_n & \jmx_n \cr \jmx_n & -\imx_n}
 \pmatrix {\emx_n & \emx_n \cr \emx_n & \emx_n}
 \pmatrix {\imx_n & \jmx_n \cr \jmx_n & -\imx_n}
 = \pmatrix{2 \emx_n & \omx_n \cr \omx_n & \omx_n}.
\nonumber\end{align}
\end{description}

\par\medskip\noindent
\subsection*{Odd Dimensional Case}
\par\medskip\noindent
We now consider
$(2n+1) \times (2n+1)$
matrices,
$n \in {\mathbb N}.$
Let
\begin{align}
 \xmx_{2n+1} &= \pmatrix{\rc{\sqrt 2}\imx_n & 0_n & \rc{\sqrt 2}\jmx_n \cr 0_n^T & 1 & 0_n^T \cr
         \rc{\sqrt 2}\jmx_n & 0_n & -\rc{\sqrt 2}\imx_n}
 \in {\mathbb R}^{(2n+1) \times (2n+1)}.
\nonumber\end{align}
(The matrix
$\xmx_{2n+1}$
turns into
$\xmx_{2n}$
when its central row and column are deleted.) The matrix
$\xmx_{2n+1}$
is a symmetric involution, i.e.
$\xmx_{2n+1}^T = \xmx_{2n+1}$
and
$\xmx_{2n+1}^2 = \imx_{2n+1}.$
Conjugation with
$\xmx_{2n+1}$
gives rise to the following block representation of odd-dimensional semimagic square
matrices.
\begin{Theorem}\label{tblro}
Let
$M \in {\mathbb R}^{(2n+1) \times (2n+1)}.$
\begin{description}
\item{(a)}
The matrix
$M \in A_{2n+1}^o$
if and only if
\begin{align}
 M &= \xmx_{2n+1} \pmatrix{\omx_n & 0_n & V^T \cr 0_n^T & 0 & -\sqrt 2 (V 1_n)^T \cr
     W & -\sqrt 2 W 1_n & \omx_n} \xmx_{2n+1}
\nonumber\end{align}
with matrices
$V, W \in {\mathbb R}^{n \times n}.$
Moreover,
$M$
will have rational entries if and only if
$V$
and
$W$
have rational entries.
\item{(b)}
The matrix
$M \in B_{2n+1}^o$
if and only if
\begin{align}
 M &= \xmx_{2n+1} \pmatrix{Y & -\sqrt 2 Y 1_n & \omx_n \cr -\sqrt 2 (Y^T 1_n)^T & 2\,1_n^T Y 1_n
   & 0_n^T \cr \omx_n & 0_n & Z} \xmx_{2n+1}
\nonumber\end{align}
with matrices
$Y, Z \in {\mathbb R}^{n \times n}.$
Moreover,
$M$
will have rational entries if and only if
$Y$
and
$Z$
have rational entries.
\item{(c)}
The matrix
$\emx_{2n+1}$
satisfies
\begin{align}
 \emx_{2n+1} &= \xmx_{2n+1} \pmatrix{2 \emx_n & \sqrt 2\,1_n & \omx_n \cr \sqrt 2\,1_n^T & 1 & 0_n^T
     \cr \omx_n & 0_n & \omx_n} \xmx_{2n+1}.
\nonumber\end{align}
\end{description}
\end{Theorem}
\par\medskip\noindent
Note that there are no conditions on the matrices
$V, W, Y$
and
$Z$
in Theorem \ref{tblro}, so the block representation gives a very simple way of constructing
all odd-dimensional semimagic square matrices (with or without centre-point symmetries).
Indeed, it is evident from the theorem that the general element of
$S_{2n+1}$
will be of the form
\begin{align}
 M &= \xmx_{2n+1} \pmatrix{Y + 2 w \emx_n & \sqrt 2 (w \imx_n - Y) 1_n & V^T \cr
   \sqrt 2 ((w \imx_n - Y) 1_n)^T & w + 2\,1_n^T Y 1_n & - \sqrt 2 (V 1_n)^T \cr
   W & -\sqrt 2 W 1_n & Z} \xmx_{2n+1}
\label{eoddform}\end{align}
with arbitrary
$V, W, Y, Z \in {\mathbb R}^{n \times n}.$
Note that, in contrast to the even-dimensional case, adding the weight
$w$
to
$M$
is not equivalent to adding a weight to the matrix
$Y.$
Incidentally, the choices $n = 1$, $w = 5$, $V = (2)$, $W = (4)$ and
$Y = Z = (0)$ give the Loh Shu square.
\par\medskip\noindent
{\it Proof\/}
of Theorem \ref{tblro}.
Writing
$M$
in the form
\begin{align}
 M &= \pmatrix{A & v & C \cr w^T & x & y^T \cr B & z & D}
\nonumber\end{align}
with
$x \in {\mathbb R},$
$v, w, y, z \in {\mathbb R}^n$
and
$A, B, C, D \in {\mathbb R}^{n \times n},$
we find
\begin{align}
 \jmx_{2n+1} M \jmx_{2n+1} &= \pmatrix{\jmx_n D \jmx_n & \jmx_n z & \jmx_n B \jmx_n \cr
     y^T \jmx_n & x & w^T \jmx_n \cr \jmx_n C \jmx_n & \jmx_n v & \jmx_n A \jmx_n}
\nonumber\end{align}
\begin{description}
\item{(a)}
The weight 0 association condition of Lemma \ref{lsemi} (b) now takes the form
\begin{align}
 \omx_{2n+1} &= M + \jmx_{2n+1} M \jmx_{2n+1}
\nonumber\\
 &= \pmatrix{\jmx_n(\jmx_n A \jmx_n + D)\jmx_n & v + \jmx_n z & \jmx_n(\jmx_n C \jmx_n + B)\jmx_n \cr
     w^T + y^T \jmx_n & 2 x & (w^T + y^T \jmx_n)\jmx_n \cr
     B + \jmx_n C \jmx_n & \jmx_n(v + \jmx_n z) & D + \jmx_n A \jmx_n},
\nonumber\end{align}
so
$x = 0,$
$y = -\jmx_n w,$
$z = -\jmx_n v,$
$B = -\jmx_n C \jmx_n$
and
$D = -\jmx_n A \jmx_n.$
Thus, if M is a weight 0 associated semimagic square matrix, then
\begin{align}
% xmx_{2n+1} M xmx_{2n+1} &= rc 2 pmatrix{imx_n & 0_n & jmx_n cr 0_n^T & sqrt 2 & 0_n^T cr
%    jmx_n & 0 & -imx_n}
% pmatrix{A & v & C cr w^T & 0 & -w^T jmx_n cr -jmx_n C jmx_n & - jmx_n v & -jmx_n A jmx_n}
% pmatrix{imx_n & 0_n & jmx_n cr 0_n^T & sqrt 2 & 0_n^T cr jmx_n & 0 & -imx_n}
 \xmx_{2n+1} M \xmx_{2n+1} &= \xmx_{2n+1} \pmatrix{A & v & C \cr w^T & 0 & -w^T \jmx_n \cr
   -\jmx_n C \jmx_n & - \jmx_n v & -\jmx_n A \jmx_n} \xmx_{2n+1}
\nonumber\\
 &= \pmatrix{\omx_n & 0_n & A \jmx_n - C \cr 0_n^T & 0 & \sqrt 2 w^T \jmx_n \cr
    \jmx_n A + \jmx_n C \jmx_n & \sqrt 2 \jmx_n v & \omx_n}.
\nonumber\end{align}
By Lemma \ref{lsemi} (a),
$M \in S_{2n+1}^o$
if and only if
$1_{2n+1}$
is an eigenvector with eigenvalue 0 of
$M$
and of
$M^T.$
Since
\begin{align}
 \xmx_{2n+1} 1_{2n+1} &= \pmatrix{\sqrt 2\,1_n \cr 1 \cr 0_n},
\label{eoddone}\end{align}
we see that
\begin{align}
 0_{2n+1} &= M 1_{2n+1} = M \xmx_{2n+1} \xmx_{2n+1} 1_{2n+1}
\nonumber\end{align}
if and only if
\begin{align}
 0_n &= W \sqrt 2\,1_n + \sqrt 2 \jmx_n v,
\nonumber\end{align}
where we set
$W = \jmx_n A + \jmx_n C \jmx_n;$
this gives
$\jmx_n v = - W 1_n.$
Similarly, setting
$V = (A \jmx_n - C)^T,$
\begin{align}
 0_{2n+1} &= M^T 1_{2n+1} = M^T \xmx_{2n+1} \xmx_{2n+1} 1_{2n+1}
\nonumber\end{align}
if and only if
\begin{align}
 0_n &= V \sqrt 2\,1_n + \sqrt 2 \jmx_n w,
\nonumber\end{align}
and hence
$w^T \jmx_n = -(V 1_n)^T.$
\hfill\break
Conversely, given matrices
$V, W \in {\mathbb R}^{n \times n},$
we take
$A = \rc 2(\jmx_n W + V^T \jmx_n),$
$C = \rc 2(\jmx_n W \jmx_n - V^T)$
and further
$B = -\jmx_n C \jmx_n,$
$D = -\jmx_n A \jmx_n,$
$x = 0,$
$v = -\jmx_n W 1_n,$
$w = -\jmx_n V 1_n,$
$y = V 1_n$
and
$z = W 1_n$
to construct a weight 0 associated semimagic square matrix.
\hfill\break
From these formulae, it is evident that
$M$
has rational entries if and only if
$V, W$
have rational entries.
\item{(b)}
The balance condition of Lemma \ref{lsemi} (c) takes the form
\begin{align}
 \omx_{2n+1} &= M - \jmx_n M \jmx_n
\nonumber\\
 &= \pmatrix{\jmx_n(\jmx_n A \jmx_n - D)\jmx_n & v - \jmx_n z & \jmx_n(\jmx_n C \jmx_n - B)\jmx_n \cr
   w^T - y^T \jmx_n & 0 & (y^T \jmx_n - w^T)\jmx_n \cr
   B - \jmx_n C \jmx_n & \jmx_n(\jmx_n z - v) & D - \jmx_n A \jmx_n},
\nonumber\end{align}
so
$z = \jmx_n v,$
$y = \jmx_n w,$
$B = \jmx_n C \jmx_n$
and
$D = \jmx_n A \jmx_n;$
there is no condition on
$x.$
Then
\begin{align}
 \xmx_{2n+1} M \xmx_{2n+1} &= \pmatrix{A + C \jmx_n & \sqrt 2 v & \omx_n \cr
   \sqrt 2 w^T & x & 0_n^T \cr \omx_n & 0_n & \jmx_n A \jmx_n - \jmx_n C}
\nonumber\end{align}
and hence, using Lemma \ref{lsemi} (a) and (\ref{eoddone}),
$0_{2n+1} = M 1_{2n+1}$
if and only if
\begin{align}
 0_{n+1} &= \pmatrix{\sqrt 2 ((A + C \jmx_n) 1_n + v) \cr 2 w^T 1_n + x},
\nonumber\end{align}
giving
$x = -2 w^T 1_n$
and
$v = -Y 1_n,$
where we set
$Y = A + C \jmx_n.$
Similarly,
$0_{2n+1} = M^T 1_{2n+1}$
if and only if
\begin{align}
 0_{n+1} &= \pmatrix{\sqrt 2 (Y^T 1_n + w) \cr 2 v^T 1_n + x},
\nonumber\end{align}
so
$x = -2 v^T 1_n$
and
$w = - Y^T 1_n.$
This determines
$v$
and
$w$
and gives two conditions on
$x,$
which turn out to be each equivalent to
$x = 2\,1_n^T Y 1_n.$
\hfill\break
For the converse, we take
$A = \rc 2 (Y + \jmx_n Z \jmx_n),$
$C = \rc 2 (Y \jmx_n - \jmx_n Z),$
and further
$B = \jmx_n C \jmx_n,$
$D = \jmx_n A \jmx_n,$
$v = -Y 1_n,$
$w = -Y^T 1_n,$
$x = 2\, 1_n^T Y 1_n,$
$y = - \jmx_n Y^T 1_n$
and
$z = -\jmx_n Y 1_n$
to construct a weight 0 balanced semimagic square matrix.
\hfill\break
From these formulae it is evident that
$M \in {\mathbb Q}^{(2n+1) \times (2n+1)}$
if and only if
$Y, Z \in {\mathbb Q}^{n \times n}.$
\item{(c)}
This is a straightforward calculation.
\phantom{.}\hfill$\qed$\par\noindent
%.lp
%[This is far too long and repetitive!]
\end{description}
%\par\bigskip\noindent{\bf {[sdihed]}. Dihedral symmetries}\nobreak
\par\medskip\noindent
The unique decomposition of a semimagic square matrix into a weightless associated,
a weightless balanced semimagic square matrix and a multiple of
$\emx_n$
(where
$n$
is the dimension of the matrix)
can easily be read off its block representation. Associated and balanced semimagic
square matrices are clearly identified by the presence and position of null blocks in
their block representation.
Moreover, the block representation can be used to characterise other matrix symmetries
in a convenient manner.
For example, in the study of semimagic or magic squares, it is usually the relative
arrangement of numbers in the square which is the object of interest, not their fixed
positioning in a matrix; therefore matrices which differ only by rotation or reflection
will be identified with one another. In other words, the equivalence classes in
$S_n$
with respect to the dihedral group for the
$n \times n$
square
will be considered. In the block representation, the action of the dihedral group
translates into transposition and sign inversion of the right or bottom half of the
block matrix.
Thus we have the following result.
\begin{Theorem}\label{tdihed}
Let
$M \in S_n$
and
\begin{align}
 M &= \xmx_n \pmatrix{\tilde Y & \tilde V^T \cr \tilde W & Z} \xmx_n
\label{eblrg}\end{align}
its block representation.
Then the dihedral equivalence class of
$M$
is
\begin{align}
 &\left\{\xmx_n \pmatrix{\tilde Y & s \tilde V^T \cr t \tilde W & s t Z} \xmx_n,
  \xmx_n \pmatrix{\tilde Y^T & s \tilde W^T \cr t \tilde V & s t Z^T} \xmx_n \mid s, t \in \{1, -1\}
 \right\}.
\nonumber\end{align}
\end{Theorem}
\par\medskip\noindent
{\it Remark.\/}
If
$n$
is even, then
$\tilde Y, \tilde V, \tilde W$
are the matrices
$Y, V, W$
of Theorem \ref{tblre}; if
$n$
is odd, then they include parts of the centre row and column in the block representation
of Theorem \ref{tblro}.
(In particular,
$\tilde V$
and
$\tilde W$
will then no longer be square matrices.)
In either case,
$Z$
is the matrix denoted by the same letter in Theorems \ref{tblre}, \ref{tblro}.
\par\medskip\noindent
{\it Proof.\/}
The dihedral group is generated by the two operations of reflection along the diagonal
and reflection along the horizontal centreline.
In matrix terms, these correspond to taking the transpose of the matrix and to left
multiplication with the matrix
$\jmx_n,$
respectively.
Transposition carries over directly to the block representation, as seen in Lemma
\ref{liso} (b).
Since
$\xmx_n \jmx_n M \xmx_n = \xmx_n \jmx_n \xmx_n \xmx_n M \xmx_n,$
left multiplication of
$M$
with
$\jmx_n$
translates into left multiplication of its block representation with the matrix
\begin{align}
 \xmx_n \jmx_n \xmx_n &= \pmatrix{\imx_k & \omx_k \cr \omx_k & -\imx_k}
\nonumber\end{align}
in the even-dimensional case
$n = 2k,$
and
\begin{align}
 \xmx_n \jmx_n \xmx_n &= \pmatrix{\imx_{k+1} & \omx_{k+1, k} \cr \omx_{k, k+1} & - \imx_k}
\nonumber\end{align}
in the odd-dimensional case
$n = 2k+1;$
here
$\omx_{k,l}$
denotes the null matrix with
$k$
rows and
$l$
columns.
Hence
\begin{align}
 \jmx_n M &= \xmx_n \pmatrix{\tilde Y & \tilde V^T \cr -\tilde W & -Z} \xmx_n.
\label{eJM}\end{align}
%Clearly this matrix product inverts the sign of the bottom
% k
%rows in the block representation.
%.pp
Similarly, reflection of
$M$
along the vertical centreline corresponds to multiplication with
$\jmx_n$
on the right, and hence to inverting the sign of the rightmost
$k$
columns in the block representation. (This operation is of course derived from the
group generators as
$(\jmx_n M^T)^T = M \jmx_n.)$
\phantom{.}\hfill$\qed$\par\noindent
\section{Magic Square Matrices}
\par\medskip\noindent
For a semimagic square matrix to be magic, its two diagonals also need to add up to
the row and column sum. For associated semimagic square matrices, this is always the
case. However, for balanced semimagic squares, which, since
$\emx_n$
is magic, we can assume without loss of generality to have weight 0, this gives an
additional condition.
\par
The sum of the diagonal entries of the matrix
$M$
is equal to its trace,
$\mathop{\rm tr}\nolimits M;$
the sum of the entries on the second diagonal is the trace of the matrix after reflection along the
horizontal (or vertical) centreline, i.e. equal to
$\mathop{\rm tr}\nolimits (\jmx_n M).$
As the trace of a product of matrices is invariant under cyclic permutations of its
factors and hence the trace of a semimagic square matrix is equal to the trace of
its block representation, we see that for a weight 0 balanced semimagic square matrix
$M,$
using the representation (\ref{eblrg}),
\begin{align}
 \mathop{\rm tr}\nolimits M &= \mathop{\rm tr}\nolimits (\xmx_n \pmatrix{\tilde Y & \omx \cr \omx & Z} \xmx_n)
 = \mathop{\rm tr}\nolimits \tilde Y + \mathop{\rm tr}\nolimits Z,
\nonumber\end{align}
and by (\ref{eJM})
\begin{align}
 \mathop{\rm tr}\nolimits (\jmx_n M) &= \mathop{\rm tr}\nolimits (\xmx_n \pmatrix{\tilde Y & \omx \cr \omx & -Z} \xmx_n)
 = \mathop{\rm tr}\nolimits \tilde Y - \mathop{\rm tr}\nolimits Z.
\nonumber\end{align}
As both must vanish for the matrix to be magic, this means that the traces of
$\tilde Y$
and of
$Z$
must separately be 0.
(As all diagonal entries of the block representation of a weight 0 associated semimagic
square matrix are 0, it is evident that such matrices are always magic.)
This gives rise to the following statement.
\begin{Theorem}\label{tmagic}
\begin{description}
\item{(a)}
A
$2n \times 2n$
matrix
$M$
is a balanced magic square matrix of weight
$w$
if and only if
\begin{align}
 M &= \xmx_{2n} \pmatrix{Y & \omx_n \cr \omx_n & Z} \xmx_{2n},
\nonumber\end{align}
where
$Y$
is an
$n \times n$
semimagic square matrix of weight
$2 w$
and
$\mathop{\rm tr}\nolimits Y = 2 n w,$
and
$Z$
is any
$n \times n$
matrix with
$\mathop{\rm tr}\nolimits Z = 0.$
\item{(b)}
A
$(2n+1) \times (2n+1)$
matrix
$M$
is a balanced magic square matrix of weight
$w$
if and only if
\begin{align}
 M &= \xmx_{2n+1} \pmatrix{\tilde Y & \omx_{n+1, n} \cr \omx_{n, n+1} & Z} \xmx_{2n+1},
\nonumber\end{align}
where
\begin{align}
 \tilde Y &= \pmatrix{Y + 2 w \emx_n & \sqrt 2 (w \imx_n - Y) 1_n \cr
      \sqrt 2 ((w \imx_n - Y) 1_n)^T & w + 2\,1_n^T Y 1_n},
% M &= xmx_{2n+1} pmatrix{Y + 2 w emx_n & sqrt 2 (w imx_n - Y) 1_n & omx_n cr
%   sqrt 2 ((w imx_n - Y) 1_n)^T & w + 2\,1_n^T Y 1_n & 0_n^T cr
%   omx_n & 0_n & Z} xmx_{2n+1},
\nonumber\end{align}
%where
$Y, Z$
are any
$n \times n$
matrices such that
$\mathop{\rm tr}\nolimits Y = -2\,1_n^T Y 1_n$
(corresponding to
$\mathop{\rm tr}\nolimits \tilde Y = (2n+1) w)$
and
$\mathop{\rm tr}\nolimits Z = 0.$
\end{description}
\end{Theorem}
\par\medskip\noindent
{\it Remarks.\/}
1. The condition on the matrix
$Y$
in Theorem \ref{tmagic} (a) almost makes it a magic square matrix, but not quite. For example, the
matrix
\begin{align}
 Y &= \pmatrix{1 & 2 & -3 \cr -6 & 1 & 5 \cr 5 & -3 & -2}
\nonumber\end{align}
is a weight 0 semimagic square matrix with vanishing trace, but its second diagonal
does not add up to 0. Nevertheless, this matrix is a suitable upper block
$Y$
for the block representation of a
$6 \times 6$
balanced magic square matrix, choosing any
$3 \times 3$
matrix
$Z$
of vanishing trace for the lower block.
\par
2. Bearing in mind that
$1_n^T Y 1_n$
is just the sum of all entries of
$Y,$
the condition on the matrix
$Y$
in Theorem \ref{tmagic} (b) means that the off-diagonal entries of
$Y$
sum to
$-\frac 3 2$
times its trace. For example,
\begin{align}
 Y &= \pmatrix{1 & -2 \cr -1 & 1}
\nonumber\end{align}
has this property, and picking a matrix
$Z$
with vanishing trace, we obtain a balanced magic square matrix by the calculation
\begin{align}
 \xmx_5 \pmatrix{1 & -2 & \sqrt 2 & 0 & 0 \cr -1 & 1 & 0 & 0 & 0 \cr 0 & \sqrt 2 & -2 & 0 & 0 \cr
    0 & 0 & 0 & 1 & 3 \cr 0 & 0 & 0 & 2 & -1} \xmx_5
 &= \pmatrix{0 & 0 & 1 & -2 & 1 \cr 1 & 1 & 0 & 0 & -2 \cr 0 & 1 & -2 & 1 & 0 \cr
    -2 & 0 & 0 & 1 & 1 \cr 1 & -2 & 1 & 0 & 0}.
\nonumber\end{align}
\par
3. Theorem \ref{tmagic} gives a simple method of constructing balanced magic square
matrices. General magic square matrices can be obtained by adding any associated
semimagic (and hence magic) square matrix, or equivalently by filling in the
off-diagonal blocks in the block representation
with matrices satisfying the conditions in Theorem \ref{tblre} (a)
or \ref{tblro} (a).
\par\medskip\noindent
{\bf Definition.}
We call an
$n \times n$
magic square matrix
{\it trivial\/}
if it is a multiple of the matrix
$\emx_n.$
\par\medskip\noindent
Any associated semimagic square matrix is magic; by (\ref{ealg}), its square will be a balanced semimagic
matrix. One may
wonder whether it is again magic. For the case of even-dimensional
associated magic square matrices with rank 1 blocks (and, w.l.o.g., weight 0), the following
{\it alternative \/}
holds.
\par\medskip\noindent
%[Do we need to review the algebra of squares?]
\begin{Theorem}\label{talt}
Consider the
(non-trivial, weight 0)
$2n \times 2n$
associated magic square
\begin{align}
 M &= \xmx_{2n} \pmatrix{\omx_n & V^T \cr W & \omx_n} \xmx_{2n}
\nonumber\end{align}
with rank 1 blocks
$V, W.$
Then exactly one of the following statements is true.
\begin{description}
\item{(a)}
$M^2$
has rank 2 and is not magic, and none of its powers are magic;
\item{(b)}
$M^2$
has rank 0 or 1 and is magic and nilpotent, in fact
$M^4 = \omx_{2n}.$
\end{description}
\end{Theorem}
\par\medskip\noindent
{\it Proof. \/}
There are vectors
$u, v, x, y \in {\mathbb R}^n\setminus\{0_n\}$
such that
$V = u v^T,$
$W = x y^T.$
The block representation of
$M^2$
is
\begin{align}
 M^2 &= \xmx_{2n} \pmatrix{V^T W & \omx_n \cr \omx_n & W V^T} \xmx_{2n}
\label{emsq}\end{align}
with blocks
\begin{align}
 V^T W &= v u^T x y^T, \qquad W V^T = x y^T v u^T.
\nonumber\end{align}
If
$u^T x = 0$
and
$y^T v \neq 0$,
then
$V^T W = \omx_n$
and
$W V^T$
is a rank 1 matrix; it only has eigenvalue 0.
Indeed, any eigenvector must be in the range of the matrix, i.e. a multiple of
$x,$
which is in the null space of the matrix by hypothesis.
Since the trace of a matrix is equal to the sum of its eigenvalues, repeated
according to algebraic multiplicity, it follows that both
$V^T W$
and
$W V^T$
have vanishing trace, and by Theorem \ref{tmagic} this implies that
$M^2$
is magic.
Furthermore,
\begin{align}
 W V^T W V^T &= x (y^T v) (u^T x) (y^T v) u^T = \omx_n,
\nonumber\end{align}
so the matrix is nilpotent.
\par
The case where
$y^T v = 0$
and
$u^T x \neq 0$
is analogous.
If both
$u^T x = 0$
and
$y^T v = 0,$
then
$M^2 = \omx_{2n}.$
\par
If both
$u^T x \neq 0$
and
$y^T v \neq 0,$
then both
$V^T W$
and
$W V^T$
have a non-zero eigenvalue with eigenvector
$v \neq 0_n, x \neq 0_n,$
respectively, as
\begin{align}
 V^T W v = v (u^T x) (y^T v), &\qquad W V^T x = x (y^T v) (u^T x);
\nonumber\end{align}
in fact they both have the
{\it same\/}
eigenvalue
$(u^T x) (y^T v) \neq 0.$
Hence each block in (\ref{emsq}) has rank 1 and a non-zero trace, so
$M^2$
has rank 2 and, by Theorem \ref{tmagic}, is not magic.
For any
$N \in {\mathbb N},$
its $N$th power has block representation
\begin{align}
 M^{2N} &= \xmx_{2n} \pmatrix{(V^T W)^N & \omx_n \cr \omx_n & (W V^T)^N} \xmx_{2n},
\nonumber\end{align}
and as both blocks have the non-zero eigenvalue
$((u^T x) (v^T y))^N,$
%the matrix
$M^{2N}$
is not magic.
\phantom{.}\hfill$\qed$\par\noindent
\par\medskip\noindent
In case (b) of Theorem \ref{talt}, the matrix may have rank 1 but only eigenvalue 0.
This means that its Jordan normal form will have a single off-diagonal 1 and zeros
otherwise, so the geometric multiplicity of the eigenvalue 0 of
$M$
will be
$2n-1.$
\par
It is an open question whether Theorem \ref{talt} generalises to higher-rank matrices.
It is generally true that if
$W V^T = \omx_n,$
then
$V^T W$
has no non-zero eigenvalues, so
$M^2$
will then be magic and nilpotent; and similarly if
$V^T W = \omx_n.$
\par
We have, however, the following general result, which shows that non-trivial magic
has limited power.
\begin{Theorem}\label{tnomagicpower}
Let
$M$
be a balanced magic square matrix of size
$m \times m,$
where
$m = 2 n$
or
$m = 2n - 1.$
Then there is
$N \in \{1, \dots, n\}$
such that
$M^N$
is either trivial or not magic.
\end{Theorem}
\par\medskip\noindent
{\it Remark. \/}
The loss of magic may be temporary, as the magic property may recur when higher powers
are involved. For example, the balanced magic square matrix
\begin{align}
 \pmatrix{0 & -1 & 1 & 0 \cr 1 & 0 & 0 & -1 \cr -1 & 0 & 0 & 1 \cr 0 & 1 & -1 & 0}
 &= \xmx_4 \pmatrix{0 & 0 & 0 & 0 \cr 0 & 0 & 0 & 0 \cr 0 & 0 & 0 & 2 \cr 0 & 0 & -2 & 0} \xmx_4
\nonumber\end{align}
has non-magic even and magic odd powers.
On the other hand, triviality, i.e. the property of being a multiple of
$\emx,$
obviously persists to all higher powers.
\par\medskip\noindent
For the proof of Theorem \ref{tnomagicpower}, we require the following Vandermonde-type lemma.
\begin{Lemma}\label{lpoly}
If the numbers
$\lambda_1, \dots, \lambda_n \in {\mathbb C}$
satisfy
\begin{align}
 \sum_{j=1}^n \lambda_j^k &= 0 \qquad (k \in \{1, \dots, n\}),
\label{eapz}\end{align}
then
$\lambda_1 = \lambda_2 = \dots = \lambda_n = 0.$
\end{Lemma}
\par\medskip\noindent
{\it Proof \/}
by induction. For
$n=1$
the statement is trivial. Assume
$n \in {\mathbb N}$
is such that the statement is true for
$n-1.$
Then, considering
the polynomial
\begin{align}
 p(x) &= \prod_{j=1}^n (x - \lambda_j)
 = \sum_{k=0}^n \alpha_k\,x^k,
\nonumber\end{align}
we see that
\begin{align}
 0 &= \sum_{j=1}^n p(\lambda_j)
 = \sum_{j=1}^n \sum_{k=0}^n \alpha_k\,\lambda_j^k
 = \sum_{k=0}^n \alpha_k \sum_{j=1}^n \lambda_j^k
 = n \alpha_0
\nonumber\end{align}
by $(\ref{eapz})$, so
$\alpha_0 = 0.$
Therefore
$0$
is a root of
$p,$
so one of the
$\lambda_j$
vanishes.
By suitable renumbering we can assume without loss of generality that
$\lambda_n = 0;$
then
$\lambda_1, \dots, \lambda_{n-1}$
satisfy
the conditions $(\ref{eapz})$ for
$n-1.$
By induction,
they all vanish.
\phantom{.}\hfill$\qed$\par\noindent
\par\medskip\noindent
{\it Remark.\/}
The lemma and its proof generalise to the case where $(\ref{eapz})$
is replaced with
\begin{align}
 \sum_{j=1}^n c_j\,\lambda_j^k &= 0
 \qquad (k \in \{1, \dots, n\})
\nonumber\end{align}
with non-zero coefficients
$c_j \in {\mathbb C} \setminus\{0\}.$
\par\medskip\noindent
{\it Proof\/}
of Theorem \ref{tnomagicpower}.
We can assume without loss of generality that
$M$
has weight 0.
By Theorem \ref{tblre} (b) or Theorem \ref{tblro} (b),
\begin{align}
 M &= \xmx_m \pmatrix{\tilde Y & \omx \cr \omx & Z} \xmx_m,
\nonumber\end{align}
where
$\tilde Y$
is an
$n \times n$
matrix and
$Z$
is an
$n \times n$
or an
$(n-1) \times (n-1)$
matrix, depending on whether
$m$
is even or odd.
If we assume that the powers
$M, M^2, \dots, M^n$
are all magic, then
$\tilde Y, \tilde Y^2, \dots, \tilde Y^n$
and
$Z, Z^2, \dots, Z^n$
all have trace 0, by Theorem \ref{tmagic}.
As the eigenvalues of
$\tilde Y^k$
are just the $k$th powers of the eigenvalues of
$\tilde Y,$
this means that the eigenvalues of
$\tilde Y$
satisfy the hypothesis of Lemma \ref{lpoly}
and thus are all zero, and likewise the eigenvalues of
$Z$
all vanish.
Consequently, the Jordan normal forms of
$\tilde Y$
and of
$Z$
have all zero entries except possibly some entries 1 above the diagonal, and hence are
nilpotent; in particular,
$\tilde Y^n = \omx,$
$Z^n = \omx,$
and hence
$M^n = \omx_m.$
\phantom{.}\hfill$\qed$\par\noindent
\par\medskip\noindent
Since the square of an associated magic square is a balanced semimagic square, Theorem
\ref{tnomagicpower} has the following immediate consequence.
\begin{Corollary}
Let
$M$
be an associated magic square matrix of size
$2n \times 2n$
or
$(2n-1) \times (2n-1).$
Then there is
$N \in \{1, \dots, n\}$
such that
$M^{2N}$
is either trivial or not magic.
\end{Corollary}
\par\medskip\noindent
Of course, all odd powers of an associated magic square matrix are again associated
(semi)magic square matrices and hence, in particular, magic.
\par
This settles the matter for
associated
and for
balanced
magic square matrices; however, it does not solve the general case. For example, taking
the block form of even-dimensional matrices,
\begin{align}
 \pmatrix{Y & V^T \cr W & Z}^2 &= \pmatrix{Y^2 + V^T W & Y V^T + V^T Z \cr W Y + Z W & W V^T + Z^2},
\nonumber\end{align}
and the conditions
$0 = \mathop{\rm tr}\nolimits Y^2 + \mathop{\rm tr}\nolimits V^T W,$
$0 = \mathop{\rm tr}\nolimits Z^2 + \mathop{\rm tr}\nolimits W V^T$
do not easily relate to the properties of these matrices separately.
However, we can apply the idea of the proof of Theorem \ref{tnomagicpower} directly to
the
{\it whole\/}
matrix, using the fact that for a weight 0 magic square matrix it is
necessary, though not sufficient,
that its trace vanish. This gives the following general statement, with a less tight
upper bound on the number
$N.$
\begin{Theorem}\label{tnomagicpowerever}
Let
$M$
be a magic square. Then there is some positive integer
$N,$
not greater than the dimension of
$M,$
such that
$M^N$
is either trivial or not magic.
\end{Theorem}

\section{Ranks, Quasi-Inverses and Spectral Properties}\nobreak
In Theorem \ref{tevac}, we have seen that adding a non-zero weight increases the
rank of a weightless even-dimensional associated magic square matrix by 1.
However, this can never lead to full rank. Indeed, since the block component
matrices
$V, W$
in the block representation
\begin{align}
 &\pmatrix{2 w \emx_n & V^T \cr W & \omx_n}
\nonumber\end{align}
have row sum 0 and therefore rank
$\le n-1,$
the total rank of the matrix is at most
$2n-2$
if the weight is 0, and
$2n-1$
otherwise.
In particular, an even-dimensional associated magic square matrix is never
regular and does not have an inverse.
\par
An
$n \times n$
semimagic square matrix with weight zero will always have the vector
$1_n$
in its null space and therefore cannot be regular; however we can define
a (left or right)
{\it quasi-inverse\/}
to be a matrix which multiplies the square (to the left or right) to give
the weightless part of the (semimagic) unit matrix,
\begin{align}
 \umx_n := \imx_n - \rc n\,\emx_n.
\nonumber\end{align}
\par
An even-dimensional associated magic square matrix, weighted or not, will not
have a left nor a right quasi-inverse, as can be seen from the block
representation
\begin{align}
 \umx_{2n} &= \xmx_{2n} \pmatrix{\umx_n & \omx_n \cr \omx_n & \imx_n} \xmx_{2n};
\label{equnit}\end{align}
the upper right-hand block of a right quasi-inverse of the block representation would have to be a right
inverse of
$W,$
while the lower left-hand block of a left quasi-inverse would have to be a left
inverse of
$V,$
and clearly neither is possible.
\par
However, in the odd-dimensional case, the block representation is (see Theorem \ref{tblro})
\begin{align}
 &\pmatrix{2 w \emx_n & w \sqrt 2\,1_n & V^T \cr
         w \sqrt 2\,1_n^T & w & -\sqrt 2 (V 1_n)^T \cr
         W & -\sqrt 2 W 1_n & \omx_n},
\nonumber\end{align}
where the matrices
$V, W$
may have full rank
$n.$
Hence the maximal rank is
$2n$
if the weight is 0 and
$2n + 1$
otherwise; in the latter case, the matrix has full rank and therefore an
inverse.
\par
For a rank
$2n,$
$(2n+1) \times (2n+1)$
associated magic square matrix with weight 0, a right quasi-inverse can always
be constructed, bearing in mind that
$V, W$
have full rank
$n$
and the matrix
$\imx_n + 2 \emx_n$
is regular, as
$-\rc 2$
is not an eigenvalue of
$\emx_n.$
Indeed,
\begin{align}
 \pmatrix{0 & 0 & V^T \cr 0 & 0 & -\sqrt{2}(V 1_n)^T \cr W & -\sqrt{2} W 1_n & 0}
 \pmatrix{0 & 0 & {V'}^T \cr 0 & 0 & -\sqrt{2}({V'} 1_n)^T \cr W' & -\sqrt{2} W' 1_n & 0}
 &= \xmx_{2n+1} \umx_{2n+1} \xmx_{2n+1},
\nonumber\end{align}
taking
$W' := (V^T)^{-1}(1 - \frac{2 \emx_n}{2n+1}),$
${V'}^T := (1+2 \emx_n)^{-1} W^{-1}.$
Note that the quasi-inverse is again the block structure matrix of a weight 0
associated magic square matrix.
Here
\begin{align}
 \xmx_{2n+1} \umx_{2n+1} \xmx_{2n+1}
 &= \pmatrix{1- \frac{2 \emx_n}{2n+1} & -\sqrt{2}\frac{1_n}{2n+1} & \omx_n \cr
            -\sqrt{2}\frac{1_n^T}{2n+1} & 1 - \rc{2n+1} & 0_n^T \cr
            \omx_n & 0_n & \imx_n}.
\nonumber\end{align}
In summary, we have the following statement.
\begin{Theorem}\label{tApseu}
\begin{description}
\item{(a)}
Even-dimensional associated magic square matrices, with or without weight,
never have full rank, nor a quasi-inverse.
\item{(b)}
Odd-dimensional associated magic square matrices may have full rank if weighted;
if the weight is 0, then the maximal rank is 1 less than the dimension, and in
this case left and right quasi-inverses exist.
\end{description}
\end{Theorem}
\par\medskip\noindent
Turning to the case of (weightless) balanced semimagic square matrices, we note
that in the block representation for the odd-dimensional case,
\begin{align}
 &\pmatrix{Y & -\sqrt 2 Y 1_n & 0 \cr -\sqrt 2 1_n^T Y & 2\,1_n^T Y 1_n & 0 \cr 0 & 0 & Z},
\nonumber\end{align}
the top left
$(n+1) \times (n+1)$
matrix
$\tilde Y$
has maximal rank
$n$
(when
$Y$
has rank
$n$
), because the $n+1$st row is a linear combination of the first
$n$
rows. Therefore the maximal rank of the balanced semimagic square, achieved
if both
$Y, Z$
have full rank
$n,$
is
$2n,$
and there is no inverse. However, there is always a quasi-inverse in this case;
with
$Y' := (1+2 \emx_n)^{-1} Y^{-1} (I_n - \frac{2}{2n+1} \emx_n),$
\begin{align}
 \pmatrix{Y & -\sqrt 2 Y 1_n \cr -\sqrt 2\,1_n^T Y & 2\,1_n^T Y 1_n}
 \pmatrix{Y' & -\sqrt 2 Y' 1_n \cr -\sqrt 2\,1_n^T Y' & 2\,1_n^T Y' 1_n}
 &= \pmatrix{\imx_n - \frac{2}{2n+1} \emx_n & \frac{-\sqrt 2}{2n+1} 1_n \cr
            \frac{-\sqrt 2}{2n+1} 1_n^T  & 1 - \rc{2n+1}},
\nonumber\end{align}
and completing a
$(2n+1) \times (2n+1)$
matrix with
$Z^{-1}$
in the lower right corner, we obtain the block representation of a right
quasi-inverse which is again a weightless balanced semimagic square matrix.
\par
The block representation of a
$2n$-dimensional weightless balanced semimagic square matrix is
\begin{align}
 &\pmatrix{Y & 0 \cr 0 & Z},
\nonumber\end{align}
where
$Y$
is a weightless semimagic square matrix. Hence the maximal possible rank is
$2n-1,$
and there will not be an inverse. Regarding a quasi-inverse, we see from
(\ref{equnit}) that
$Z$
must be invertible and
$Y$
must have a quasi-inverse.
If
$n$
is odd and
$Y$
has maximal rank
$n-1$
and is either associated or balanced, then there exists a quasi-inverse by
the above considerations.
If
$n$
is even and
$Y$
is associated, then no quasi-inverse exists.
If
$n$
is even and
$Y$
is balanced, then we can apply these considerations recursively to the block structure of
$Y.$
\par
Unfortunately, the case of mixed type, i.e. of a general weightless semimagic
square matrix, seems rather more difficult to analyse;
this case will generally occur when applying
the above reduction to an even-dimensional balanced semimagic square matrix.
Also, the case of a balanced semimagic square matrix with weight seems rather difficult, as it is not obvious
whether the addition of a weight always raises the rank by 1.

\begin{Theorem}\label{tsqaddweight}
Let
\begin{align}
 M_0 &= \xmx_{2n} \pmatrix{\omx_n & V^T \cr W & \omx_n} \xmx_{2n}
\nonumber\end{align}
be a weightless associated magic square matrix, and for
$w \in {\mathbb R} \setminus\{0\},$
let
$M_w := M_0 + w \emx_{2n}.$
Then
$M_w^2$
has the same non-zero eigenvalues (including multiplicities) as
$M_0^2,$
with $w$-independent eigenvectors,
and the additional eigenvalue
$2 n^2 w^2.$
Moreover,
$\mathop{\rm rk} M_w^2 = \mathop{\rm rk} M_0^2 + 1.$
\end{Theorem}
\par\medskip\noindent
{\it Proof.\/}
Since
$M_0 \emx_{2n} = \omx_{2n} = \emx_{2n} M_0,$
\begin{align}
 M_w^2 &= M_0^2 + w\,(M_0 \emx_{2n} + \emx_{2n} M_0) + w^2 \emx_{2n}^2
 = M_0^2 + 2n w^2 \emx_{2n}.
\nonumber\end{align}
Now assume that
$\lambda \neq 0$
is an eigenvalue of the block representation of
$M_0^2$
(which, as
$\xmx_{2n}$
is unitary,
has the same eigenvalues as
$M_0^2$
),
so there is a vector
\begin{align}
 w &= \pmatrix{w_1 \cr w_2} \in {\mathbb C}^{2n} \setminus\{0\}
\nonumber\end{align}
such that
\begin{align}
 \pmatrix{ V^T W & \omx_n \cr \omx_n & W V^T } \pmatrix{w_1 \cr w_2}
 &= \lambda \pmatrix{w_1 \cr w_2}.
\nonumber\end{align}
Then
$V^T W w_1 = \lambda w_1,$
so
$w_1 \in \mathop{\rm ran} V^T,$
i.e.
$w_1$
is a linear combination of columns of
$V^T,$
and hence
$1_n^T w_1 = 0.$
Consequently,
$\emx_{n} w_1 = 0_n.$
Therefore, considering the block representation of
$M_w^2,$
\begin{align}
 \pmatrix{ V^T W + 4n w^2 \emx_n & \omx_n \cr \omx_n & W V^T } \pmatrix{w_1 \cr w_2}
 &= \pmatrix{ \lambda w_1 + 0_n \cr \lambda w_2} = \lambda \pmatrix{w_1 \cr w_2},
\nonumber\end{align}
which shows that
$\lambda$
is still an eigenvalue of
$M_w^2.$
Moreover, using the fact that
$W 1_n = 0_n,$
we see that
\begin{align}
 \pmatrix{ V^T W + 4n w^2 \emx_n & \omx_n \cr \omx_n & W V^T } \pmatrix{ 1_n \cr 0_n }
 &= 4 n^2 w^2 \pmatrix{ 1_n \cr 0_n },
\nonumber\end{align}
so
$4 n^2 w^2$
is a further eigenvalue of
$M_w^2.$
\par
Conversely, if
$\lambda$
is an eigenvalue of
$M_w^2,$
so
\begin{align}
 \pmatrix{ V^T W + 4n w^2 \emx_n & \omx_n \cr \omx_n & W V^T } \pmatrix{ w_1 \cr w_2}
 &= \lambda \pmatrix{ w_1 \cr w_2},
\nonumber\end{align}
then we can write
$w_1 = \alpha 1_n + z,$
where
$z \in {\mathbb R}^n$
with
$z^T 1_n = 0,$
and find
\begin{align}
 (V^T W + 4 n w^2 \emx_n) w_1 &= V^T W z + 4 n^2 w^2 \alpha 1_n = \lambda z + \lambda \alpha 1_n.
\nonumber\end{align}
Since
$V^T W z$
is orthogonal to
$1_n,$
this implies
\begin{align}
 V^T W z &= \lambda z, \qquad \lambda \alpha = 4 n^2 w^2 \alpha,
\nonumber\end{align}
so
$\lambda = 4 n^2 w^2$
unless
$\alpha = 0,$
and
$\lambda$
is an eigenvalue of
$M_0^2$
unless
$z = w_2 = 1_n.$
\par
For the statement about ranks, note that
$V^T W$
maps the orthogonal complement
$\{1_n\}^\bot \subset {\mathbb R}^n$
into itself, so adding a non-zero multiple of
$\emx_n$
adds one dimension, parallel to
$1_n,$
to the range of the matrix.
\phantom{.}\hfill$\qed$\par\noindent
\par\medskip\noindent
We remark that the following analogous statement holds for the eigenvalues of
the associated magic square matrix itself.
\begin{Theorem}\label{tevac}
Let
$M_0$
and
$M_w$
be as in Theorem \ref{tsqaddweight}. Then
$M_w$
has the same non-zero eigenvalues as
$M_0$
and an additional simple eigenvalue
$2 n w.$
The eigenvectors are independent of
$w.$
Moreover,
$\mathop{\rm rk} M_w = \mathop{\rm rk} M_0 + 1.$
\end{Theorem}
\par\medskip\noindent
{\it Proof. \/}
Assume
\begin{align}
 \pmatrix{0 & V^T \cr W & 0} \pmatrix{w_1 \cr w_2} &= \lambda \pmatrix{w_1 \cr w_2}
\nonumber\end{align}
with non-zero
$\lambda;$
then
$V^T w_2 = \lambda w_1$
%and
% W w_1 = lambda w_2.
%In particular,
implies that
$1_n^T w_1 = 0.$
The remainder of the proof is completely analogous to that
of Theorem \ref{tsqaddweight}.
\phantom{.}\hfill$\qed$\par\noindent
\par\medskip\noindent
Hence we can omit the weight of the associated magic square w.l.o.g.\ in the following.
The weightless squared block matrix
\begin{align}
 \pmatrix{V^T W & \omx_n \cr \omx_n & W V^T}&,
\nonumber\end{align}
is obviously the direct sum of the matrices
$V^T W$
and
$W V^T,$
with separate quadratic forms --- the
$\xi$
part and the
$\eta$
part of the quadratic form (\ref{eqf}) ---
and correspondingly the quadratic form generated by
$M^2$
splits in a natural way into two quadratic forms
which are generated by
$B A$
and
$A B,$
respectively.
These forms will be expressible as a sum of squares (possibly with irrational coefficients)
if both
$V^T W$
and
$W V^T$
are diagonalisable, in particular if they are symmetric.
This motivates the following definition.
\par\medskip\noindent
{\bf Definition.}
We call the associated magic square matrix
$M$
{\it parasymmetric\/}
if its square
$M^2$
is a symmetric matrix.
\par\medskip\noindent
In terms of the block representation, parasymmetry can be characterised as follows if the
two constituent blocks of
$M$
have rank 1.
\begin{Lemma}\label{lparasym}
Let
$V, W \in {\mathbb R}^{n \times n}$
have row sum
$0$
and rank 1, and consider the associated magic square matrix
\begin{align}
 M &= \xmx_{2n} \pmatrix{\omx_n & V^T \cr W & \omx_n} \xmx_{2n}.
\nonumber\end{align}
If
$M^2 \neq \omx_{2n},$
then
$M$
is parasymmetric if and only if
$W$
is a multiple of
$V.$
\end{Lemma}
\par\medskip\noindent
{\it Proof.\/}
We can write
$V = u v^T,$
$W = x y^T$
with non-null vectors
$u, v, x, y \in {\mathbb R}^n \setminus \{0_n\};$
then
$V^T = v u^T$
and
$W^T = y x^T.$
Now if
$V^T W$
is symmetric, then
\begin{align}
 (u^T x) v y^T = V^T W = W^T V = (x^T u) y v^T&,
\nonumber\end{align}
so either
$u^T x = 0$
(and note that
$u^T x = x^T u$
) or
$v y^T = y v^T;$
but in the latter case we see, multiplying by
$y$
on the right, that
\begin{align}
 v (y^T y) = y (v^T y),
\nonumber\end{align}
so
$v$
and
$y$
are linearly dependent.
Similarly, if
$W V^T$
is symmetric, then either
$v^T y = 0$
or
$u$
and
$x$
are linearly dependent.
\par
Now, if it should happen that
$u^T x = 0,$
then we must also have
$v^T y = 0$
(since otherwise, by the above,
$u, x$
will be simultaneously orthogonal and linearly dependent, which is impossible as both
are non-null vectors), and vice versa; and in this situation
$V^T W = \omx_n = W V^T,$
which would imply
$M^2 = \omx_{2n}.$
\par
Hence, we find that
$x, y$
are multiples of
$u, v,$
respectively, so
$W$
is a multiple of
$V.$
\par
The converse statement is obvious.
\phantom{.}\hfill$\qed$\par\noindent
\par\medskip\noindent
More generally, a matrix is diagonalisable by conjugation with an orthogonal matrix if it
commutes with its transpose, i.e. if it is
{\it normal.\/}
Obviously, any symmetric matrix is normal.
If
$M$
is normal, then it is easy to see that
$M^2$
is normal, too; the converse is not so clear. In analogy to parasymmetry, we call the matrix
$M$
{\it paranormal\/}
if
$M^2$
is normal. However, this turns out to be no more general than parasymmetry as far as
associated magic square matrices with rank 1 blocks are concerned, as the following result
shows.
\begin{Lemma}\label{lparanorm}
Let
$M \in {\mathbb R}^{2n \times 2n}$
be a weightless associated magic square matrix with rank 1 blocks
$V, W$
such that
$M^2 \neq \omx_{2n}.$
If
$M$
is paranormal, then it is parasymmetric.
\end{Lemma}
\par\medskip\noindent
{\it Proof.\/}
Let
$V, W$
and
$u,v,x,y$
be as in the proof of Lemma \ref{lparasym}. Then the paranormality means that
\begin{align}
 V^T W W^T V &= W^T V V^T W,
\nonumber\\
 W V^T V W^T &= V W^T W V^T;
\nonumber\end{align}
in terms of the generating vectors, this gives the two identities
\begin{align}
 v (u^T x)(y^T y)(x^T u) v^T &= y (x^T u)(v^T v)(u^T x) y^T,
\nonumber\\
 x (y^T v)(u^T u)(v^T y) x^T &= u (v^T y)(x^T x)(y^T v) u^T,
\nonumber\end{align}
and by the same reasoning as above, this implies that
$u^T x = 0$
or
$v, y$
are linearly dependent, and that
$v^T y = 0$
or
$x, u$
are linearly dependent. As before, the cases of orthogonality can only occur together
and then give a trivial
$M^2;$
it follows that
$V, W$
are linearly dependent and hence that the matrix is parasymmetric.
\phantom{.}\hfill$\qed$\par\noindent
\par\medskip\noindent
The following observation on the eigenvalues of a squared weightless
associated magic square matrix with rank 1 blocks follows from the proof
of Theorem \ref{talt}.
\begin{Theorem}\label{tinteigen}
If $M$ is the weightless associated magic square matrix
\begin{equation}
\label{erankopo}
 M = \xmx_{2n} \pmatrix{\omx_{n} & v u^T \\ x y^T & \omx_n} \xmx_{2n}
\end{equation}
with rank 1 block components
$V = u v^T,$
$W = x y^T,$
then
$M^2$
has eigenvalues
$0$
and
$(u^T x)(y^T v) = \mathop{\rm Tr} V^T W.$
In particular, if
$V, W$
have integer entries, then the eigenvalues of
$M^2$
are integers.
\end{Theorem}
\par\medskip\noindent
{\it Remarks.\/}
1.
If the eigenvalue
$(u^T x)(y^T v)$
is non-zero, then
its (algebraic and geometric) multiplicity
will be 2, as it will be an eigenvalue of both the upper left and the lower
right blocks in the block representation of
$M^2,$
with eigenvectors
$\pmatrix{v \cr 0_n}$
and
$\pmatrix{0_n \cr x},$
respectively.
However, when exactly one of the products
$V^T W, W V^T$
vanishes,
then its geometric multiplicity will only be 1.
This can happen in the non-parasymmetric case, as the non-orthogonality of
generating vectors for non-trivial magic squares, found in the proof of Lemma \ref{lparasym},
need not hold in this case.
Note that in this situation, this eigenvalue will indeed be 0, so the matrix
$M^2$
will have eigenvalue 0 only with algebraic multiplicity
$2n$
and geometric multiplicity
$2n - 1.$
For example, consider
\begin{align}
 V &= u v^T = \pmatrix{1 \cr -3 \cr -5 \cr 7} (1, -1, -1, 1),
\nonumber\\
 W &= x y^T = 8 \pmatrix{1 \cr -3 \cr -5 \cr 7} (1, -1, 1, -1);
\nonumber\end{align}
here
$x$
is a multiple of
$u,$
but
$y^T v = 0,$
so
$W V^T = 0.$
The resulting weightless associated magic square matrix
\begin{align}
 M &= \rc 2 \pmatrix{
  63 & -61 &  53 & -55 & -57 &  59 & -51 &  49 \cr
 -47 &  45 & -37 &  39 &  41 & -43 &  35 & -33 \cr
 -31 &  29 & -21 &  23 &  25 & -27 &  19 & -17 \cr
  15 & -13 & 5 &  -7 &  -9 &  11 &  -3 & 1 \cr
  -1 & 3 & -11 & 9 & 7 &  -5 &  13 & -15 \cr
  17 & -19 &  27 & -25 & -23 &  21 & -29 &  31 \cr
  33 & -35 &  43 & -41 & -39 &  37 & -45 &  47 \cr
 -49 &  51 & -59 &  57 &  55 & -53 &  61 & -63 \cr
 }
\nonumber\end{align}
has rank 2, but its square
$M^2$
only has rank 1.
\par
2. If
$(u^T x)(y^T v) \neq 0,$
then corresponding linearly independent (right) eigenvectors of
$M^2$
will be
\begin{align}
 \sqrt 2 \xmx_{2n} \pmatrix{v \cr 0_n} &= \pmatrix{v \cr \jmx_n v}, \qquad
 -\sqrt 2 \xmx_{2n} \pmatrix{0_n \cr x} = \pmatrix{-\jmx_n x \cr x};
\label{ereigen}\end{align}
the first of these is even, the other odd under reflection (i.e. multiplication
with
$\jmx_{2n}).$
These eigenvectors are orthogonal for structural reasons, reflecting the fact that
they belong to different direct summands in the block representation of $M^2$.
Clearly, if
$v, x$
have integer entries, then so do these eigenvectors.

We note that there are also the left eigenvectors (i.e. eigenvectors of
$(M^2)^T$) given by
\begin{align}
 \pmatrix{y \cr \jmx_n y}, &\pmatrix{-\jmx_n u \cr u};
\label{eleigen}\end{align}
again, these eigenvectors are structurally orthogonal.
\par
3. In the parasymmetric case
$y = v,$
$x = k u,$
where
$u, v \in {\mathbb R}^n \setminus\{0\},$
the matrix
$M^2$
always has a non-zero eigenvalue
$k (u^T u)(v^T v).$
\par
4.
In the situation of Theorem \ref{tinteigen} with
$u^T x \neq 0,$
$y^T v \neq 0,$
the matrix
$M$
has the two simple non-zero eigenvalues
$\sqrt{(u^T x) (y^T v)}, -\sqrt{(u^T x) (y^T v)}.$
Indeed, any eigenvector of the block representation
\begin{align}
 \pmatrix{\omx_n & v u^T \cr x y^T & \omx_n}&
\nonumber\end{align}
for non-zero eigenvalue
$\lambda$
can easily be seen to be of the form
$\pmatrix{\alpha v \cr \beta x}$
with
$\alpha, \beta \neq 0,$
and hence
$\lambda^2 = (u^T x) (y^T v).$
As the trace of the matrix vanishes and
$0$
is the only other eigenvalue, both signs of the square root occur.
The coefficients for the eigenvectors can be chosen as
$\alpha = \lambda,$
$\beta = y^T v,$
so if
$u, v, x, y$
are integer vectors and
$\lambda$
is an integer, then there are integer eigenvectors as well.
In the parasymmetric case with integer vectors
$u, v$
and integer parasymmetry factor
$k,$
the eigenvalues will be integers if and only if the square-free part of
$k$
is equal to the square-free part of the product
$(v^T v)(u^T u).$

\begin{Lemma}\label{minpoly}
Let
$n \ge 2$,
and let
$M \in {\mathbb R}^{2n \times 2n}$ be a weightless associated magic square matrix with $n\times n$ rank 1 block components
\[
M= {\xmx}_{2n}\pmatrix{\omx_n & v u^T \cr x y^T & \omx_n}{\xmx}_{2n},
\]
where
$u, v, x, y \in \mathbb{R}^n$
such that
$\lambda := (u^T x)(y^T v) \neq 0$.
Then the minimal polynomial of $M$ is
\[
 x^3 - \lambda x.
\]
\end{Lemma}
\par\medskip\noindent
{\it Remark.\/}
In the case $n = 1$, the minimal polynomial is
$ x^2 - \lambda.$
\par\medskip\noindent
{\it Proof.\/}
By straightforward calculation,
\begin{align*}
M^3 &= {\xmx}_{2n}\pmatrix{\omx_n & v u^T \cr x y^T & \omx_n} \pmatrix{\omx_n & v u^T \cr x y^T & \omx_n} \pmatrix{\omx_n & v u^T \cr x y^T & \omx_n}{\xmx}_{2n}
\\
 &= {\xmx}_{2n}\pmatrix{\omx_n & v (u^T x)(y^T v) u^T \cr x (y^T v)(u^T x) y^T & \omx_n}{\xmx}_{2n},
\\
 &= (u^T x)(y^T v){\xmx}_{2n}\pmatrix{\omx_n & v u^T \cr x y^T & \omx_n}{\xmx}_{2n} = \lambda M.
\end{align*}
To see that this is indeed the minimal polynomial for $M$, we note that
the block representation of $M$ is not a multiple of the unit matrix, which rules out
a linear polynomial, and that when inserted in a quadratic polynomial
$x^2 + a x + b$,
\[
\pmatrix{v(u^T x) y^T & \omx_n \cr \omx_n & x(y^T v) u^T}
+ a \pmatrix{\omx_n & v u^T \cr x y^T & \omx_n} + b \pmatrix{\imx_n & \omx_n \cr \omx_n & \imx_n}
= \omx_{2n}
\]
implies %$a = 0$ and
$v(u^T x) y^T = -b \imx_n = x(y^T v) u^T$, which is impossible, as $\imx_n$ has
rank $n > 1$.
\phantom{.}\hfill$\qed$\par\noindent

\medskip\noindent
To conclude this section, we show a construction yielding a two-sided, regular eigenvector matrix for $M^2$, where $M$ is the rank $1 + 1$ associated magic square matrix as defined in (\ref{erankopo}).
Here `two-sided' means that the columns of the matrix are right eigenvectors of
$M^2$ while its rows are left eigenvectors of $M^2$.
We begin by considering the two right eigenvectors (\ref{ereigen}) of $M^2$ corresponding to the non-zero eigenvalue $\lambda = (u^T x)(y^T v)$.
These eigenvectors are placed side by side to form a
$2n \times 2$
matrix
\begin{align}
 P_1 &= \pmatrix{B_1 \cr A_1 \cr C_1} = \pmatrix{v & -\jmx_n x \\ \jmx_n v & x},
\nonumber\end{align}
where
 $A_1 = \pmatrix{v_n & -x_1 \cr v_n & x_1}$
and
$B_1$ and $C_1$ are $(n-1) \times 2$ matrices such that
$C_1 = \jmx_{n-1} B_1 \sigma_3$
(with
$\sigma_3 = \pmatrix{1 & 0 \cr 0 & -1}$).
We make the assumption that
$v_n, x_1 \neq 0$,
so that the matrix
$A_1$
is regular.
The construction below can be generalised to the case where $A_1$ is any regular
matrix composed of two rows of $P_1$ (and indeed one can always find two linearly
independent rows of $P_1$ because its columns are linearly independent), but we
take the centre rows in the following for simplicity.

Now define
\begin{align}
 \tilde P_1 &= - P_1 A_1^{-1} = \pmatrix{-B_1 A_1^{-1} \cr -\imx_2 \cr -C_1 A_1^{-1}}.
\label{ePonetil}\end{align}
Similarly, starting from the matrix of left eigenvectors (\ref{eleigen}), we set
\begin{align}
 P_2 &= \pmatrix{B_2 \cr A_2 \cr C_2} = \pmatrix{y & -\jmx_n u \\ \jmx_n y & u},
\nonumber\end{align}
and assuming $y_n, u_1 \neq 0$, so that
 $A_2 = \pmatrix{y_n & -u_1 \cr y_n & u_1}$
is regular, we define
\begin{align}
 \tilde P_2 &= - P_2 A_2^{-1} = \pmatrix{-B_2 A_2^{-1} \cr -\imx_2 \cr -C_2 A_2^{-1}}.
\label{ePtwotil}\end{align}
The columns of $\tilde P_1$ and $\tilde P_2$ will still be linearly independent
eigenvectors, for eigenvalue $\lambda$, of $M^2$ and of $(M^2)^T$, respectively,
but in general they will no longer be orthogonal unless $v_n^2 = x_1^2$ and $y_n^2 = u_1^2$.
However, due to the special structure of our chosen matrices $A_j$, they will
have the symmetry that the second column is the reversal of the first.
Indeed, $\jmx_2 A_j \sigma_3 = A_j$, so $\sigma_3 A_j^{-1} \jmx_2 = A_j^{-1}$
$(j \in \{1, 2\})$, and it follows that
\begin{align*}
 \jmx_{2n} \tilde P_j \jmx_2 &= \pmatrix{-\jmx_{n-1} C_j A_j^{-1} \cr -\jmx_2 \cr
     -\jmx_{n-1} B_j A_j^{-1}} \jmx_2
 = \pmatrix{-B_j \sigma_3 A_j^{-1} \jmx_2 \cr -\imx_2 \cr -C_j \sigma_3 A_j^{-1}
      \jmx_2}
  = \tilde P_j
  \qquad (j \in \{1, 2\}),
\end{align*}
which means that swapping the columns of $\tilde P_j$ is tantamount to turning
them upside down.

The matrices $\tilde P_1$ and $\tilde P_2$ have the following remarkable
connection with the matrix $M$.
%
%\begin{Lemma}
%Let $u, v, x, y \in \mathbb{R}^n$ such that $(u^T x)(y^T v) \neq 0$ and
%$u_1, v_n, x_1, y_n \neq 0$. Let $M$ be the matrix (\ref{erankopo}), and let
%$\tilde P_1$ and $\tilde P_2$ be defined as in (\ref{ePonetil}), (\ref{ePtwotil}).
%Then
%\[
%(\omx_{2n,n-1} \mid \tilde P_1 \mid\omx_{2n,n-1}) M=-M
%\]
%and
%\[
%(\omx_{2n,n-1} \mid \tilde P_2 \mid \omx_{2n,n-1}) M^T = - M^T,
%\]
%where $(\cdot \mid \cdot \mid \cdot)$ just means that the three matrices are
%juxtaposed to form one $2n \times 2n$ matrix.
%\end{Lemma}
%
%\smallskip\noindent
%{\it Proof.}
Using the notation $\left(\cdot \mid \cdot \mid \cdot\right)$ to express that the three
matrices are juxtaposed to form one $2n \times 2n$ matrix,
we calculate
\begin{align}
 (\omx_{2n,n-1} \mid \tilde P_1 \mid\ &\omx_{2n,n-1}) M
 = - P_1 (\omx_{2,n-1} \mid A_1^{-1} \mid \omx_{2,n-1}) \xmx_{2n} \pmatrix{\omx_n & v u^T \cr
   x y^T & \omx_n} \xmx_{2n}
\nonumber\\
 &= - \rc 2\,P_1 \left( \omx_{2,n-1} \mid A_1^{-1} \xmx_2 \mid \omx_{2,n-1} \right)
  \pmatrix{\omx_n & v u^T \cr x y^T & \omx_n} \pmatrix{\imx_n & \jmx_n \cr \jmx_n & -\imx_n}
\nonumber\\
 &= - \rc 2\,P_1 \left( A_1^{-1} \pmatrix{x_1 \cr -x_1} y^T \mid A_1^{-1} \pmatrix{v_n \cr v_n} u^T \right) \pmatrix{\imx_n & \jmx_n \cr \jmx_n & -\imx_n}
\nonumber\\
 &= - \rc 2 \pmatrix{\imx_n & \jmx_n \cr \jmx_n & -\imx_n} \pmatrix{v & 0_n \cr 0_n & -x}
  \pmatrix{0_n^T & u^T \cr -y^T & 0_n^T} \pmatrix{\imx_n & \jmx_n \cr \jmx_n & -\imx_n}
\nonumber\\
 &= - \xmx_{2n} \pmatrix{ \omx_n & v u^T \cr x y^T & \omx_n} \xmx_{2n} = - M,
\nonumber\end{align}
and similarly
$(\omx_{2n,n-1} \mid \tilde P_2 \mid \omx_{2n,n-1}) M^T = - M^T.$
%\phantom{.}\hfill$\qed$\par\noindent

\begin{Theorem}\label{teigen}
Let $u, v, x, y \in \mathbb{R}^n$ such that $\lambda := (u^T x)(y^T v) \neq 0$ and
$u_1, v_n, x_1, y_n \neq 0$. Let $M$ be the matrix (\ref{erankopo}), and let
$\tilde P_1$ and $\tilde P_2$ be defined as in (\ref{ePonetil}), (\ref{ePtwotil}).
Then
\begin{align}
 P &= \imx_{2n} + (\omx_{2n,n-1} \mid \tilde P_1 \mid \omx_{2n,n-1}) +
  \pmatrix{\omx_{n-1,2n} \cr \tilde P_2^T \cr \omx_{n-1,2n}}
\nonumber\end{align}
is a two-sided eigenvector matrix for
$M^2,$
so that
$M^2 P = P \mathop{\rm diag}(0_{n-1}, \lambda 1_2, 0_{n-1})$
and
$P M^2 = \mathop{\rm diag}(0_{n-1}, \lambda 1_2, 0_{n-1}) P.$
Moreover,
$P$
has the inverse
\begin{align}
 P^{-1} &= \mathop{\rm diag}(1_{n-1}, 0_2, 1_{n-1}) - \frac{M^2} \lambda.
\nonumber\end{align}
\end{Theorem}

\smallskip\noindent
{\it Proof.}
Using the second of the above identities and the fact that the columns of
$\tilde P_1$ are eigenvectors of $M^2$, we find
\begin{align}
 M^2 P &= M^2 + M^2 (\omx_{2n,n-1} \mid \tilde P_1 \mid \omx_{2n,n-1}) + M M
  \pmatrix{\omx_{n-1,2n} \cr \tilde P_2^T \cr \omx_{n-1,2n}}
\nonumber\\
 &= M^2 + \lambda (\omx_{2n,n-1} \mid \tilde P_1 \mid \omx_{2n,n-1}) + M (-M);
\nonumber\end{align}
on the other hand, since the central
$2 \times 2$
part of
$\tilde P_2$
is equal to
$-\imx_2,$
we have
\begin{align}
 P \mathop{\rm diag}(0_{n-1}, \lambda 1_2, 0_{n-1})
 &= \mathop{\rm diag}(0_{n-1}, \lambda 1_2, 0_{n-1}) + (\omx_{2n,n-1} \mid \tilde P_1 \mid \omx_{2n,n-1}) \mathop{\rm diag}(0_{n-1}, \lambda 1_2, 0_{n-1})
\nonumber\\
 &\qquad + \pmatrix{\omx_{n-1,2n} \cr \tilde P_2^T \cr \omx_{n-1,2n}} \mathop{\rm diag}(0_{n-1}, \lambda 1_2, 0_{n-1})
\nonumber\\
 = \mathop{\rm diag}(0_{n-1}, &\lambda 1_2, 0_{n-1}) + \lambda (\omx_{2n,n-1} \mid \tilde P_1 \mid \omx_{2n,n-1}) - \mathop{\rm diag}(0_{n-1}, \lambda 1_2, 0_{n-1}).
\nonumber\end{align}
The relation for $ P M^2$ follows by a pair of completely analogous calculations.

To verify the inverse relation
\begin{align}
 P^{-1} &= \mathop{\rm diag}(1_{n-1}, 0_2, 1_{n-1}) - \frac{M^2} \lambda,
\nonumber\end{align}
we note that
\begin{align}
 &(\mathop{\rm diag}(1_{n-1}, 0_2, 1_{n-1}) - \frac{M^2} \lambda) P
 = \mathop{\rm diag}(1_{n-1}, 0_2, 1_{n-1}) P - P \mathop{\rm diag}(0_{n-1},1_2,0_{n-1})
\nonumber\\
 &= \mathop{\rm diag}(1_{n-1}, 0_2, 1_{n-1}) + \mathop{\rm diag}(1_{n-1}, 0_2, 1_{n-1})\,(\omx_{2n,n-1} \mid \tilde P_1 \mid \omx_{2n,n-1})
\nonumber\\
 &\qquad + \mathop{\rm diag}(1_{n-1}, 0_2, 1_{n-1}) \pmatrix{\omx_{n-1,2n} \cr \tilde P_2^T \cr \omx_{n-1,2n}}
 - \mathop{\rm diag}(0_{n-1},1_2,0_{n-1})
\nonumber\\
 &\qquad - (\omx_{2n,n-1} \mid \tilde P_1 \mid \omx_{2n,n-1})\,\mathop{\rm diag}(0_{n-1},1_2,0_{n-1}) - \pmatrix{\omx_{n-1,2n} \cr \tilde P_2^T \cr \omx_{n-1,2n}}\mathop{\rm diag}(0_{n-1},1_2,0_{n-1})
\nonumber\\
 &= \mathop{\rm diag}(1_{n-1}, 0_2, 1_{n-1}) + \Big( (\omx_{2n,n-1} \mid \tilde P_1 \mid \omx_{2n,n-1}) + \mathop{\rm diag}(0_{n-1},1_2,0_{n-1}) \Big) + \omx_{2n}
\nonumber\\
 &\qquad - \mathop{\rm diag}(0_{n-1},1_2,0_{n-1}) - (\omx_{2n,n-1} \mid \tilde P_1 \mid \omx_{2n,n-1}) + \mathop{\rm diag}(0_{n-1},1_2,0_{n-1})
\nonumber\\
 &= \mathop{\rm diag}(1_{2n}) = \imx_{2n}.
\nonumber\end{align}
The opposite product follows similarly.
\phantom{.}\hfill$\qed$\par\noindent

\begin{remark}
In the
{\it parasymmetric\/}
case, we have
$y = v, x = k u,$
which, following the above construction, gives rise to the eigenvector matrices
\begin{align}
 P_1 &= \pmatrix{y & -\jmx_n k u \cr \jmx_n y & k u} = P_2 \pmatrix{1 & 0 \cr 0 & k}.
\nonumber\end{align}
However, in this situation the vector
$\pmatrix{-\jmx_n u \cr u}$
will be an eigenvector just as well as
$\pmatrix{-\jmx_n k u \cr k u},$
so we can begin with this vector and take
$P_1 = P_2$
instead of the above.
Hence, in this instance,
$P = \imx_{2n} + (\omx_{n-1} \mid \tilde P_1 \mid \omx_{n-1}) + (\omx_{n-1} \mid \tilde P_2 \mid \omx_{n-1})^T$
will be a
{\it symmetric\/}
matrix.
\end{remark}
\section{Quadratic Forms from Squares of Associated Magic Squares}\nobreak
In this section we focus on the case of
$2n \times 2n$
associated magic squares, establishing a connection between their block representation vectors and certain types of quadratic forms.
Here, our aim is just to establish a link between the associated magic square matrices with rank $1+1$ block representations and quadratic forms. A deeper examination of how the vector structures and quadratic forms interrelate over both the field of rationals and the ring of integers will be left to later further investigation.

We start from the observation that the block representation of the associated magic
square
\begin{align}
 M &= \xmx_{2n} \pmatrix{2 w \emx_n & V^T \cr W & \omx_n} \xmx_{2n}
\nonumber\end{align}
has a natural decomposition into 3 parts, as
\begin{align}
 \pmatrix{2 w \emx_n  & V^T \cr W & \omx_n}
 &= 2 w \pmatrix{ \emx_n  & \omx_n \cr \omx_n & \omx_n} + \pmatrix{\omx_n & V^T \cr \omx_n & \omx_n} + \pmatrix{\omx_n & \omx_n \cr W & \omx_n} =: 2 w e + a + b,
\nonumber\end{align}
where
$a^2 = \omx_n,$
$b^2 = \omx_n,$
$e^2 = n e,$
$a e = \omx_n,$
$e b = \omx_n,$
$e a = \omx_n,$
$b e = \omx_n$
(the last two identities following from the fact that the rows of
$V$
and
$W$
sum to 0), and
\begin{align}
 a b = \pmatrix{V^T W & \omx_n \cr \omx_n & \omx_n}, \qquad b a &= \pmatrix{\omx_n & \omx_n \cr \omx_n & W V^T}.
\nonumber\end{align}
This corresponds to a splitting of the magic square
\begin{align}
 M &= w \emx_{2n} + A + B,
\nonumber\end{align}
where
\begin{align}
 A &= \xmx_{2n} a \xmx_{2n} = \rc 2 \pmatrix{V^T J & - V^T \cr J V^T J & -J V^T},
\qquad
 B = \xmx_{2n} b \xmx_{2n} = \rc 2 \pmatrix{J W & J W J \cr -W & -W J}.
\nonumber\end{align}
Now if we consider the (balanced semimagic square) matrix
\begin{align}
 M^2 = \xmx_{2n} \pmatrix{4 n w^2 \emx_n + V^T W & \omx_n \cr \omx_n & W V^T} \xmx_{2n},
\nonumber\end{align}
its block representation generates the quadratic form on
${\mathbb R}^{2n}$
\begin{align}
 \pmatrix{\xi \cr \eta}^T \xmx_{2n} M^2 \xmx_{2n} \pmatrix{\xi \cr \eta}
 &= 4 n w^2 \xi^T \emx_n \xi + \xi^T V^T W \xi + \eta^T W V^T \eta
 \qquad (\xi, \eta \in {\mathbb R}^n).
\label{eqf}\end{align}
Splitting the vector
$\xi$
into a part parallel to
$1_n$
and a part orthogonal to
$1_n,$
\begin{align}
 \xi &= \rc n (1_n^T \xi) 1_n + (\xi - \rc n (1_n^T \xi) 1_n),
\nonumber\end{align}
we find
\begin{align}
 \emx_n \xi&= \rc n (1_n^T \xi) n 1_n + 0_n
\nonumber\end{align}
and
\begin{align}
 W \xi &= 0_n + W (\xi - \rc n (1_n^T \xi) 1_n);
\nonumber\end{align}
this means that the first term in the quadratic form (\ref{eqf})
%$(@\eqf)$
only acts non-trivially in the subspace of
${\mathbb R}^n$
spanned by
$1_n,$
the second term only in the orthogonal subspace. Hence, by effectively restricting the
quadratic form to the subspace
$1_n^T \oplus {\mathbb R}^n,$
we can fully separate off the weight and essentially reduce the study to the weightless matrix.
(This is related to the spectral separation property of the weight matrix
observed in Theorem \ref{tsqaddweight}.)
%\par\medskip\noindent
\par
Considering more specifically the quadratic form (\ref{eqf}) with weight
$w = 0$
and taking for
$M$
a parasymmetric associated magic square matrix with rank 1 blocks ((\ref{erankopo})
with $y = v$, $x = k u$), we find
\begin{align}
 \xi^T V^T W \xi + \eta^T W V^T \eta &= k (u^T u) (v^T \xi)^2 + k (v^T v) (u^T \eta)^2,
\nonumber\end{align}
so setting
\begin{align}
 \pmatrix{\xi \cr \eta} &= \xmx_{2n} X
 = \rc{\sqrt 2} \pmatrix{X_1 + \jmx_n X_2 \cr \jmx_n X_1 - X_2}
\nonumber\end{align}
for
$X = \pmatrix{X_1 \cr X_2} \in \mathbb{R}^{2n} = \mathbb{R}^n \oplus \mathbb{R}^n$,
we obtain the quadratic form
\begin{align}
 X^T M^2 X^T &= \frac k 2\,(u^T u) \left(\pmatrix{v \cr \jmx_n v}^T X \right)^2
 + \frac k 2\,(v^T v) \left(\pmatrix{-\jmx_n u \cr u}^T X \right)^2
 \qquad
 (X \in {\mathbb R}^{2n}).
\nonumber\end{align}
Introducing the reduced variables
\begin{align}
 z_1 &= \pmatrix{v \cr \jmx_n v}^T w, \qquad
 z_2 = \pmatrix{-\jmx_n u \cr u}^T X,
\nonumber\end{align}
this gives the reduced quadratic form
\begin{align}
 q_1(z_1, z_2) = \frac k 2\,\left((u^T u)\,z_1^2 + (v^T v)\,z_2^2 \right).
\label{eqone}\end{align}

Alternatively we can represent the quadratic form in terms of the spectral
decomposition of $M^2$, still focussing on the case of a rank $1 + 1$ matrix
$M$ (\ref{erankopo}), but not necessarily parasymmetric.
We begin with any regular (right) eigenvector matrix $P\in \mathbb{R}^{2n \times 2n}$
which has two columns,
$b_1, b_2$
(not necessarily the first two columns)
which are eigenvectors of
$M^2$
for the non-zero eigenvalue
$\lambda,$
while the other columns
$b_3, \dots, b_{2n}$
are eigenvectors for eigenvalue 0.
\par
Then, with the transformation
$X = P \alpha = \sum\limits_{j=1}^{2n} \alpha_j b_j$
(where we sort the indices of the vector
$\alpha \in {\mathbb R}^{2n}$
in the manner corresponding to the numbering of the columns of
$P$), we find the quadratic form
\begin{align}
 q(\alpha) &= X^T M^2 X = (\alpha_1 b_1 + \alpha_2 b_2 + \sum_{j=3}^{2n}\alpha_j b_j)^T \lambda\,(\alpha_1 b_1 + \alpha_2 b_2)
\nonumber\\
 &= \lambda \left(b_1^T b_1 \alpha_1^2 + 2 b_1^T b_2 \alpha_1 \alpha_2 + b_2^T b_2 \alpha_2^2
  + \sum_{j=3}^{2n} (b_j^T b_1 \alpha_j \alpha_1 + b_j^T b_2 \alpha_j \alpha_2) \right).
\label{eqnull}\end{align}
In the parasymmetric case where
$M^2$
is symmetric,
$b_3, \dots, b_{2n}$
are orthogonal on
$b_1, b_2,$
so the quadratic form simplifies to the binary form
\begin{align}
 q_2(\alpha_1, \alpha_2) = X^T M^2 X &= \lambda (b_1^T b_1 \alpha_1^2 + 2 b_1^T b_2 \alpha_1 \alpha_2 + b_2^T b_2 \alpha_2^2);
\label{eqtwo}\end{align}
if also
$b_1$
and
$b_2$
are mutually orthogonal (as is the case for the `natural' eigenvectors (\ref{ereigen}), but not
necessarily for the transformed column eigenvectors of the matrix
$P$
constructed in Theorem \ref{teigen}), then we get the simple form
\begin{align}
 q_3(\alpha_1, \alpha_2) = X^T M^2 X &= \lambda(b_1^T b_1 \alpha_1^2 + b_2^T b_2 \alpha_2^2).
\label{eqthree}\end{align}

We now illustrate these results with a parasymmetric and a non-parasymmetric example.

\medskip\noindent
{\bf Example 1.}
Let $u^T=(11,-13,-19,21)$, $x=2u$ and $v^T=y^T=(-1,1,1,-1)$ be our vectors in $\mathbb{R}^4$, and
\[
M=\xmx_8\left (\begin{array}{cc}
\omx_4 & v u^T \\
x y^T & \omx_4
\end{array}\right )\xmx_8
\footnotesize
=\frac{1}{2}\left(
\begin{array}{cccccccc}
 -63 & 61 & 55 & -53 & -31 & 29 & 23 & -21 \\
 59 & -57 & -51 & 49 & 27 & -25 & -19 & 17 \\
 47 & -45 & -39 & 37 & 15 & -13 & -7 & 5 \\
 -43 & 41 & 35 & -33 & -11 & 9 & 3 & -1 \\
 1 & -3 & -9 & 11 & 33 & -35 & -41 & 43 \\
 -5 & 7 & 13 & -15 & -37 & 39 & 45 & -47 \\
 -17 & 19 & 25 & -27 & -49 & 51 & 57 & -59 \\
 21 & -23 & -29 & 31 & 53 & -55 & -61 & 63 \\
\end{array}
\right).
\normalsize
\]
Then $M^2$ has the block representation
\[
\footnotesize
M^2=8\xmx_8\left(
\begin{array}{cccccccc}
 273 & -273 & -273 & 273 & 0 & 0 & 0 & 0 \\
 -273 & 273 & 273 & -273 & 0 & 0 & 0 & 0 \\
 -273 & 273 & 273 & -273 & 0 & 0 & 0 & 0 \\
 273 & -273 & -273 & 273 & 0 & 0 & 0 & 0 \\
 0 & 0 & 0 & 0 & 121 & -143 & -209 & 231 \\
 0 & 0 & 0 & 0 & -143 & 169 & 247 & -273 \\
 0 & 0 & 0 & 0 & -209 & 247 & 361 & -399 \\
 0 & 0 & 0 & 0 & 231 & -273 & -399 & 441 \\
\end{array}
\right)\xmx_8.
\normalsize
\]
Formula (\ref{eqone}) immediately gives the quadratic form
\begin{equation*}
 q_1(z_1, z_2) = 1092 z_1^2 + 4 z_2^2.
\end{equation*}
The matrix $M^2$ has the eigenvectors (\ref{ereigen})
\begin{equation*}
 b_1 = \pmatrix{v \cr \jmx_n v}, \qquad b_2 = \pmatrix{-\jmx_n u \cr u}
\end{equation*}
(where we divided the vector $b_2$ by $k$, as suggested in the Remark after
Theorem \ref{teigen}) for eigenvalue
$\lambda = k(u^T u) (v^T v) = 8736$; so
$b_1^T b_1 = 2 v^T v = 8$, $b_2^T b_2 = 2 u^T u = 2184$, and these vectors are
orthogonal. This gives rise to the quadratic form (\ref{eqthree})
\begin{equation*}
 q_3(\alpha_1, \alpha_2) = 8736 \,(8 \alpha_1^2 + 2184 \alpha_2^2)
 = 17472 \,(4 \alpha_1^2 + 1092 \alpha_2^2).
\end{equation*}

Applying the construction of Theorem \ref{teigen}, we obtain the rational symmetric (left and right) eigenvector matrix
\[
\footnotesize
P=\frac{1}{11}\left(
\begin{array}{cccccccc}
 11 & 0 & 0 & -16 & 5 & 0 & 0 & 0 \\
 0 & 11 & 0 & 15 & -4 & 0 & 0 & 0 \\
 0 & 0 & 11 & 12 & -1 & 0 & 0 & 0 \\
 -16 & 15 & 12 & -11 & 0 & -1 & -4 & 5 \\
 5 & -4 & -1 & 0 & -11 & 12 & 15 & -16 \\
 0 & 0 & 0 & -1 & 12 & 11 & 0 & 0 \\
 0 & 0 & 0 & -4 & 15 & 0 & 11 & 0 \\
 0 & 0 & 0 & 5 & -16 & 0 & 0 & 11 \\
\end{array}
\right).
\normalsize
\]
Taking the middle columns of this matrix for the (non-orthogonal) eigenvectors $b_1, b_2 (= \jmx_8 b_1)$,
we now have
$b_1^T b_1 = b_2^T b_2 = \frac{788}{121}$, $b_1^T b_2 = -\frac{304}{121}$,
and formula (\ref{eqtwo}) gives the quadratic form
\begin{equation*}
 q_2(\alpha_1, \alpha_2) = 8736\,\left(\frac{788}{121} \alpha_1^2 - \frac{608}{121}
 \alpha_1 \alpha_2 + \frac{788}{121} \alpha_2^2\right)
 = \frac{34944}{121}\,(197 \alpha_1^2 - 152 \alpha_1 \alpha_2 + 197 \alpha_2^2).
\end{equation*}
We note that Theorem \ref{teigen} also gives the inverse for $P$,
\[
P^{-1}=\text{diag}(1_3,0_2,1_3)-\frac{M^2}{\lambda}=\frac{8}{\lambda}
\footnotesize
\left(
\begin{array}{cccccccc}
 735 & 336 & 273 & -252 & -21 & 0 & -63 & 84 \\
 336 & 775 & -260 & 241 & 32 & -13 & 44 & -63 \\
 273 & -260 & 871 & 208 & 65 & -52 & -13 & 0 \\
 -252 & 241 & 208 & -197 & -76 & 65 & 32 & -21 \\
 -21 & 32 & 65 & -76 & -197 & 208 & 241 & -252 \\
 0 & -13 & -52 & 65 & 208 & 871 & -260 & 273 \\
 -63 & 44 & -13 & 32 & 241 & -260 & 775 & 336 \\
 84 & -63 & 0 & -21 & -252 & 273 & 336 & 735 \\
\end{array}
\right),
\normalsize
\]
and $P^{-1}M^2 P=PM^2P^{-1}=\text{diag}(0_3,\lambda,\lambda,0_3)$.

Note that $q_1$, $q_2$ and $q_3$, as constructed above, are three different
quadratic forms. They are of course all equivalent (and indeed equivalent to the
simple circular binary $x_1^2 + x_2^2$) on $\mathbb{R}^2$.
The forms $q_1$ and $q_3$ are still equivalent on $\mathbb{Q}^2$, via the variable transformation
$z_1 = 8 \alpha_1$, $z_2 = 2184 \alpha_2$,
however, they are not equivalent on $\mathbb{Z}^2$, as their discriminants differ
(cf. \cite{gauss} \S157, \cite{dirichlet} \S56).

\medskip\noindent
In our second example, we consider the block decompositions and non-positive-definite quadratic form that results from our rational eigenvector transformation applied to a non-parasymmetric type A matrix.

\medskip\noindent
{\bf Example 2.}
Let $u^T=(10,-14,-18,22)$, $x^T=(23,-25,-39,41)$ and $v^T=y^T=(-1,1,1,-1)\in\mathbb{R}^4$, and
\[
M=\xmx_8\left (\begin{array}{cc}
\omx_4 & v u^T \\
x y^T & \omx_4
\end{array}\right )\xmx_8
\footnotesize
=\left(
\begin{array}{cccccccc}
 -63 & 59 & 55 & -51 & -31 & 27 & 23 & -19 \\
 61 & -57 & -53 & 49 & 29 & -25 & -21 & 17 \\
 47 & -43 & -39 & 35 & 15 & -11 & -7 & 3 \\
 -45 & 41 & 37 & -33 & -13 & 9 & 5 & -1 \\
 1 & -5 & -9 & 13 & 33 & -37 & -41 & 45 \\
 -3 & 7 & 11 & -15 & -35 & 39 & 43 & -47 \\
 -17 & 21 & 25 & -29 & -49 & 53 & 57 & -61 \\
 19 & -23 & -27 & 31 & 51 & -55 & -59 & 63 \\
\end{array}
\right).
\normalsize
\]
Then $M^2$ has the block representation
\[
\footnotesize
M^2=8\xmx_8\left(
\begin{array}{cccccccc}
 273 & -273 & -273 & 273 & 0 & 0 & 0 & 0 \\
 -273 & 273 & 273 & -273 & 0 & 0 & 0 & 0 \\
 -273 & 273 & 273 & -273 & 0 & 0 & 0 & 0 \\
 273 & -273 & -273 & 273 & 0 & 0 & 0 & 0 \\
 0 & 0 & 0 & 0 & 115 & -161 & -207 & 253 \\
 0 & 0 & 0 & 0 & -125 & 175 & 225 & -275 \\
 0 & 0 & 0 & 0 & -195 & 273 & 351 & -429 \\
 0 & 0 & 0 & 0 & 205 & -287 & -369 & 451 \\
\end{array}
\right)\xmx_8,
\normalsize
\]
which is clearly non-symmetric.
Here we cannot use the formulae obtained above for the parasymmetric
case, but we can construct the quadratic form from the general formula (\ref{eqnull}).
The non-zero eigenvalue is $\lambda = (u^T x)(y^T v) = 8736$, and the (right and
left) eigenvector matrix constructed as in Theorem \ref{teigen} has the form
\[
\footnotesize
P=\frac{1}{115}\left(
\begin{array}{cccccccc}
 115 & 0 & 0 & -160 & 45 & 0 & 0 & 0 \\
 0 & 115 & 0 & 155 & -40 & 0 & 0 & 0 \\
 0 & 0 & 115 & 120 & -5 & 0 & 0 & 0 \\
 -184 & 161 & 138 & -115 & 0 & -23 & -46 & 69 \\
 69 & -46 & -23 & 0 & -115 & 138 & 161 & -184 \\
 0 & 0 & 0 & -5 & 120 & 115 & 0 & 0 \\
 0 & 0 & 0 & -40 & 155 & 0 & 115 & 0 \\
 0 & 0 & 0 & 45 & -160 & 0 & 0 & 115 \\
\end{array}
\right).
\normalsize
\]
Now ordering the columns of $P$ as
$P = (b_3, b_4, b_5, b_1, b_2, b_8, b_7, b_6)$ to define eigenvectors $b_1, \dots, b_8$, we see that $b_2 = \jmx_8 b_1$ and
$b_{j+3} = \jmx_8 b_j$ $(j \in \{3, 4, 5\})$, and calculate further
$b_1^T b_1 = b_2^T b_2 = \frac{80900}{115^2}$, $b_1^T b_2 = -\frac{28000}{115^2}$
as well as $b_3^T b_1 = -b_3^T b_2 = -b_6^T b_1 = b_6^T b_2 = \frac{2760}{115^2}$,
$b_4^T b_1 = -b_4^T b_2 = -b_7^T b_1 = b_7^T b_2 = -\frac{690}{115^2}$ and
$b_5^T b_1 = -b_5^T b_2 = -b_8^T b_1 = b_8^T b_2 = -\frac{2070}{115^2}$.
Hence
\begin{align}
 \sum_{j=3}^8 (b_j^T b_1 \alpha_j \alpha_1 + b_j^T b_2 \alpha_j \alpha_2)
 &= \sum_{j=3}^5 (b_j^T b_1\,(\alpha_j - \alpha_{j+3})\,\alpha_1
   + b_j^T b_2\,(\alpha_j - \alpha_{j+3})\,\alpha_2)
\nonumber\\
 &= \sum_{j=3}^5 (b_j^T b_1)\,(\alpha_j - \alpha_{j+3})\,(\alpha_1 - \alpha_2);
\nonumber
\end{align}
in our example,
\begin{equation*}
\sum_{j=3}^5 (b_j^T b_1)\,(\alpha_j - \alpha_{j+3}) =
 \frac{6}{115}\,(4\,(\alpha_3 - \alpha_6) - (\alpha_4 - \alpha_7)
 - 3\,(\alpha_5 - \alpha_8)),
\end{equation*}
and in view of the single variable $\alpha_7$ we can see that, irrespective of
the number field or ring from which the variables $\alpha_3, \dots, \alpha_8$ are
taken, their combined effect is that of a single variable, say $\alpha_0$, from
that field or ring.

In summary, we obtain the ternary quadratic form (\ref{eqnull})
\begin{equation*}
 q(\alpha_1, \alpha_2, \alpha_0) = 8736 \left(\frac{4}{529}\,(809 \alpha_1^2 - 560 \alpha_1 \alpha_2
 + 809 \alpha_2^2) + \frac{6}{115}\,\alpha_0\,(\alpha_1 - \alpha_2)\right).
\end{equation*}

\bigskip\noindent
{\bf Acknowledgement.\/}
Sally Hill's research was supported by EPSRC DTP grant EP/L504749/1.
%\par\bigskip\noindent{\bf References.}\nobreak

\end{document}